\documentclass[12pt]{amsart}
\usepackage{amsmath,amssymb,latexsym,amsthm, amscd}
\usepackage{mathrsfs}
\usepackage[all]{xy}

\theoremstyle{definition}
\newtheorem{Def}[subsection]{Definition}%[section]
\newtheorem{example}[subsection]{Example}%[section]
\theoremstyle{plain}
\newtheorem{prop}[subsection]{Proposition}
\newtheorem{thm}[subsection]{Theorem}
\newtheorem{lem}[subsection]{Lemma}
\newtheorem{cor}[subsection]{Corollary}
\newtheorem{Qst}[subsection]{Question}

\newcommand{\mbf}{\mathbf}
\newcommand{\mrm}{\operatorname}
\newcommand{\I}{\operatorname{I}}

\newcommand{\bin}[2]{
\left[
   \begin{array}{@{}c@{}}
     #1 \\#2
   \end{array}
\right]    }
 
\title[Affine canonical basis elements]{AR-quiver approach to affine canonical basis elements}
\author{ Yiqiang Li}
\address{Department of Mathematics\\ Kansas State University\\ Manhattan, KS  66506}
\email{yqli@math.ksu.edu}
\author{Zongzhu Lin}
\address{Department of Mathematics\\ Kansas State University\\ Manhattan, KS  66506}
\email{zlin@math.ksu.edu}
\thanks{The research is partially supported by NSF grant DMS-0200673\\
The main results in the paper have been reported by 
Y. Li in the conference on Lie Algebras, Vertex Operator
Algebras and Their Applications in honor of
Robert L. Wilson and James Lepowsky, NC State University, May 17-21, 2005.}

\begin{document}

\begin{abstract}
This is the continuation of ~\cite{Li}. We describe the affine canonical basis
elements in the case when the affine quiver has arbitrary orientation. 
This generalizes the description in ~\cite{lusztig2}.
\end{abstract}

\maketitle

\section{Introduction}
Let $\mbf U^-$ be the negative part of the quantized enveloping algebra $\mbf U$
of the Kac-Moody Lie algebra
associated to a symmetric generalized Cartan matrix $C$. 
One of the milestones in the Lie theory is the discovery of Lusztig's canonical basis 
$\mbf B$ (or Kashiwara's global crystal basis ~\cite{Kashiwara}) of $\mbf U^-$. 
The canonical basis $\mbf B$ possesses many remarkable properties 
such as total positivity and integrality (see ~\cite{lusztig1} and ~\cite{lusztig3}).

There are two different approaches of defining the canonical basis $\mbf B$. 
One is algebraic and the other geometric. 
For the algebraic approach, see ~\cite{Lus2} and ~\cite{Kashiwara}. 
The geometric approach is also done in ~\cite{Lus2} when $C$ is positive definite and 
is extended to all cases in ~\cite{lusztig1}. 
It is shown in ~\cite{Grojnowski-Lusztig} that the two different approaches 
produce the same basis in $\mbf U^-$. 

Let $Q$ be a quiver such that the associated Cartan matrix is $C$. 
In the geometric approach, Lusztig studies certain perverse sheaves over the representation
space $E_{V,Q}$ of $Q$. 
An algebra $\mathcal K_Q$ is constructed and is shown to be isomorphic to $\mbf U^-$.
The simple perverse sheaves in $\mathcal K_Q$ form a basis $\mathcal B_Q$ of 
$\mathcal K_Q$. 
The canonical basis $\mbf B$ is then defined to be $\mathcal B_Q$ if one identifies
$\mathcal K_Q$ with $\mbf U^-$. In other words, the geometric approach gives a geometric
realization of $\mbf U^-$ and $\mbf B$ for each $Q$ such that the associated Cartan matrix 
is $C$.
(Note that for each $C$, 
there may be several quivers such that the associated Cartan matrices are $C$.)

It is well-known that Lie theory and the representation theory of quivers 
have very deep connections ever since Gabriel's theorem ~\cite{Gabriel} was found. 
The interaction between the two theories attracts a lot of attentions since then.
Note that $\mbf U^-$ and $\mbf B$ are objects coming purely from Lie theory. 
Although the framework of Lusztig's geometric realizations of $\mbf U^-$ and $\mbf B$
is the representation space $E_{V,Q}$ of the quiver $Q$. 
No representation theory of quivers is used explicitly. 
Nor are there clear connections between the representation theory of quivers and 
Lusztig's  geometric realizations of $\mbf U^-$ and $\mbf B$.
One may ask to what extent the representation theory of quivers can help in understanding
the quantum group $\mbf U$ and the canonical basis $\mbf B$ or $\mathcal B_Q$.
Keeping this in mind, one may ask the following natural questions. 
\begin{Qst}
\label{q1}
Characterize the elements in $\mathcal B_Q$ 
by using the representation theory of 
quivers. More precisely, describe what the supports and the corresponding
local systems of the elements in $\mathcal B_Q$ are in terms of
representations of quivers.
\end{Qst}
\begin{Qst}
\label{q2}
Recover the canonical basis via bases arising in Ringel-Hall algebras.
\end{Qst}
\begin{Qst}
\label{q3}
Show the positivity or the integrality of the canonical basis $\mbf B$
in an elementary way with the help of the representation theory of quivers, 
instead of the theory of perverse sheaves.
\end{Qst}

It seems very difficult at this moment to answer Question ~\ref{q3}. 
For an attempt to answer Question ~\ref{q2}, see ~\cite{Lus2} when $Q$
is of finite type and ~\cite{LXZ} when $Q$ is of affine type and the
references therein. 
By now there are only partial answers to Question ~\ref{q1}.

The answer to Question ~\ref{q1} is known when $Q$ is of finite type, i.e.,
when the Cartan matrix is positive definite. 
The elements in $\mathcal B_Q$ are simple perverse sheaves whose supports are
$G_V$-orbits $O$ in $E_{V,Q}$ and whose restrictions to $O$ are the constant sheaf on $O$.
The answer to Question ~\ref{q1} is known when $Q$ is a cyclic quiver. 
Elements in $\mathcal B_Q$ are simple perverse sheaves whose supports are 
the aperiodic $G_V$-orbits $O$ in $E_{V,Q}$ 
and whose restrictions to $O$ are the constant sheaf on $O$.

From now on in the introduction, 
we assume that $Q$ is an affine quiver but not a cyclic quiver.

When $Q$ is a McKay quiver, i.e., all vertices are either a sink or a source, 
the answer to Question ~\ref{q1} is given in ~\cite[Theorem 6.16]{lusztig2} by using 
the representation theory of McKay quivers. The theory of representation of
McKay quivers in \cite{lusztig2} is based on McKay's correspondence.  
In ~\cite{lusztig2}, only McKay quivers  and  cyclic quivers are considered.
Note that the construction of $\mathcal K_Q$ and $\mathcal B_Q$ applies to any quivers and
the language used in ~\cite{DR} for the representation theory of affine quivers 
works for all affine quivers, not just McKay quivers.

Given any two affine quivers $Q$ and $Q'$ such that the associated
generalized Cartan matrices coincide. Although the inverse images
$\phi^{-1}(\mathcal B_Q)$ and $(\phi')^{-1}(\mathcal B_{Q'})$ coincide
(see Theorem ~\ref{coincide}), the isomorphism classes of simple
equivariant perverse sheaves in $\mathcal B_Q$ and $\mathcal B_{Q'}$
may be quite different. This becomes obvious when $Q$ and $Q'$ are the Kronecker
quiver and the cyclic quiver with two vertices, respectively.

The goal of this paper is to give a complete answer 
to Question ~\ref{q1} for any affine quiver, by
carrying out Lusztig's arguments in ~\cite{lusztig2} to any affine quiver via 
the representation theory developed in ~\cite{DR}.

In ~\cite{Li}, such a goal  has been accomplished when $Q$ is of type $\tilde{A}$. 
Similar to [Li],
the crucial part on the way of generalizing Lusztig's argument 
is still to construct a $'$nice$'$ functor 
from the category of nilpotent representations of a cyclic quiver $\mrm C_p$ 
to the full subcategory $\mrm{Rep}(T)$ generated by a tube $T$ of period $p$
which will give an equivariant morphisms between the representation
varieties and transporting the equivariant simple perverse sheaves to
equivariant simple perverse sheaves. 
It turns out that the Hall functor in ~\cite{FMV} 
will suffice to overcome this difficulty. 
This functor gives us all the properties we need in the proof of
Lemma \ref{tube} in Section \ref{noncyclic}.

It should be mentioned that after the preprint of the paper is written, the
second author received a preprint \cite{nakajima} by Nakajima, in which he
outlines an approach in answering Question ~\ref{q1} 
when the quivers are  of type $ \tilde{D}$ and $ \tilde{E}$
by using the description of $\mathcal B_{Q}$ for McKay quivers  by Lusztig
(\cite{lusztig2} Theorem 6.16)
and then apply the reflection functors. Our argument does not depend on
Theorem 6.16 in ~\cite{lusztig2}.

{\bf Acknowledgements.} 
We wish to thank  Professor Bangming Deng and Professor  Jie Xiao 
for very helpful discussions. The second author thanks Professor Hiraku
Nakajima for discussing questions in his paper \cite{nakajima}.
The first author thanks the second author for allowing him to use the
results in this paper as part of his Ph.D. thesis.

\section{Representation theory of Affine quivers}
In this section, we give a brief review of representation theory of
affine quivers. See ~\cite{DR} and ~\cite{BGP} for more details.
\subsection{Preliminary}
\label{preliminary}
A $quiver$ is an oriented graph. 
It is a quadruple $Q=(I,\Omega, h, t)$, 
where $I$ and  $\Omega$ are two finite sets 
and $h$, $t$ are two maps from $\Omega$ to $I$ such that $h(\omega)\neq t(\omega)$ 
for any $\omega\in \Omega$. 
$I$ and $\Omega$ are called the vertex and arrow sets respectively. 
Pictorially, $t(\omega) \overset{\omega}{\to} h(\omega)$ stands 
for any arrow $\omega\in \Omega$ 
and we call $h(\omega)$ and $t(\omega)$ the head and the tail of 
the arrow $\omega$ respectively. 
An {\em affine quiver} is a quiver whose underlying graph is of type $\tilde{A}_n$,
$\tilde{D}_n$, $\tilde{E}_6$, $\tilde{E}_7$ or $\tilde{E}_8$.

{\em Throughout this paper, all quivers considered will be affine. 
We fix an algebraically closed field $K$ once and for all.}

A $representaion$ of a quiver over $K$ is a pair $(V, x)$, 
where $V=\oplus_{i\in I}V_i$ is an $I$-graded finite dimensional
$K$-vector space and 
$x=\{x_{\omega}:V_{t(\omega)}\to V_{h(\omega)}|\;\omega\in \Omega\}$ 
is a collection of linear maps.
A morphism $f: (V,x)\rightarrow (W,y)$ is a collection of linear maps 
$\{f_i:V_i\to W_i|\;i\in I\}$ such that 
$f_{h(\omega)}\;x_{\omega}=y_{\omega}\;f_{t(\omega)}$, for any $\omega\in \Omega$.   
This defines a (abelian) category, denoted by $\mrm{Rep}(Q)$.

A {\em nilpotent} representation is a representation $(V,x)$ having the property: there is 
an $N$ such that for any $\omega_1,\cdots,\omega_N$ in $\Omega$ satisfying 
$t(\omega_k)=h(\omega_{k-1})$ for any $k$, the composition of morphisms 
$x_{\omega_N}\circ \cdots\circ x_{\omega_1}$ is zero. 
Denote by $\mrm {Nil}(Q)$ the full subcategory of $\mrm{Rep}(Q)$ 
of all nilpotent representations of $Q$.  
Note that when $Q$ has no oriented cycles, every representation is
nilpotent,  
i.e., $\mrm{Nil}(Q)=\mrm{Rep}(Q)$.

Denote by $\mrm{Ind}(Q)$ the set of representatives of all pairwise
nonisomorphic  
indecomposable representations in $\mrm{Rep}(Q)$.

The Euler form $<\ ,\ >$ on $\mathbb Z[I]$ is defined by
\[<\alpha,\beta> \, =\sum_{i\in I}\alpha_i\beta_i-
\sum_{\omega\in \Omega}\alpha_{t(\omega)}\beta_{h(\omega)},\] 
for any $\alpha=\sum_i \alpha_ii$, $\beta=\sum_i \beta_ii\in \mathbb Z[I]$.
The symmetric Euler form $(,)$ on $\mathbb Z[I]$ is defined to be
\[(\alpha,\beta)=\,<\alpha,\beta>+<\beta,\alpha>,\] 
for any $\alpha=\sum_i \alpha_ii$, $\beta=\sum_i \beta_ii\in \mathbb Z[I]$. 
Given $\mbf V=(V,x)\in \mrm{Rep}(Q)$, denote by $|\mbf V|$ its dimension vector
$\sum_{i\in I}(\mrm{dim}V_i)i\in \mathbb Z[I]$.
For any $M$, $N\in \mrm{Rep}(Q)$, we have from \cite{DR}
\[<|M|,|N|>\,=\mrm{dim} \;\mrm{Hom}_Q(M,N)- \mrm{dim}\; \mrm{Ext}^1_Q(M,N).\]  

Given any $I$-graded $K$-vector space $V$, let
\[E_{V,\Omega}=\oplus_{\omega\in \Omega}\, \mrm{Hom}(V_{t(\omega)},V_{h(\omega)})
\quad \text{and}\quad
G_{V}=\prod_{i \in I} \,\mrm{GL}(V_i),\]
where $\mrm{GL}(V_i)$ is the general linear group of $V_i$ for all $i\in I$.
$G_V$ acts on $E_{V,\Omega}$ naturally, i.e.,
$g.x=y$, where $y_{\omega}=g_{h(\omega)}\,x_{\omega}\,g_{t(\omega)}^{-1}$ for all
$\omega\in \Omega$.
Note that for any $x\in E_{V,\Omega}$, $(V,x)$ is a representation of $Q$. We call
that $x$ is nilpotent if $(V,x)$ is a nilpotent representation.

\subsection{BGP reflection functors}

A vertex $i \in I$ is called a sink (resp., a source) if  
$i\in \{t(\omega),h(\omega)\}$ implies $h(\omega)=i$ (resp.,
$t(\omega)=i$) for any $\omega\in \Omega$.

For any $i \in I$,
let $\sigma_{i}Q=(I, \Omega,h', t')$ be the quiver with the orientation
$(t', h')$ defined by 
\begin{enumerate}

\item[] $t'(\omega)=t(\omega)$ and $ h'(\omega)=h(\omega)$ 
if $i\notin \{t(\omega), h(\omega)\}$;

\item[] $t'(\omega)=h(\omega)$ and $h'(\omega)=t(\omega)$ 
if $i\in \{t(\omega), h(\omega)\}$;
\end{enumerate}
for any $\omega\in \Omega$.
In other words, $\sigma_i Q$ is the quiver obtained by reversing 
the arrows in $Q$ that start or terminate at $i$.

If $i$ is a sink in $Q$, 
define a functor 
\[\Phi_{i}^+:\mrm{Rep}(Q) \to \mrm{Rep} (\sigma_{i}Q)\] 
in the following way. 
For any $(V, x) \in \mrm{Rep}(Q)$,
$\Phi_{i}^+(V,x)=(W,y)\in \mrm{Rep}(\sigma_{i}Q)$, 
where $W=\oplus_{j\in I} \, W_j$ such that
\begin{enumerate}
\item[] $W_j=V_j$, if $j\neq i$;

\item[] $W_i$ is the kernel of the linear map
$\sum_{\omega \in \Omega:\;h(\omega)=i} x_{\omega}\; : 
\oplus_{\omega\in \Omega:h(\omega)=i} V_{t(\omega)} \to V_i$;

\end{enumerate}
and $y=(y_{\omega}\; |\; \omega \in \Omega)$ such that
\begin{enumerate} 

\item[] $y_{\omega}=x_{\omega}$, if $t'(\omega)\neq i$; 

\item[] $y_{\omega}:W_i\to W_{h'(\omega)}$ is the composition of the natural maps: 
        \[W_i \to \oplus_{\omega\in \Omega:\;h(\omega)=i}V_{t(\omega)}
        =\oplus_{\omega\in \Omega:t'(\omega)=i}W_{h'(\omega)} \to
        W_{h'(\omega)},\] if $t'(\omega)=i$. 

\end{enumerate}
Note that the assignments extend to a functor.
$\Phi^+_i$ is called the BGP reflection functor with respect to the sink
$i$.

Since $Q$ has no oriented cycles, 
we can order the vertices in $I$, 
say $(i_1 ,\cdots,i_n)$ ($|I|=n$), 
in such a way that $i_r$ is a sink
in the quiver $\sigma_{i_{r-1}}\cdots \sigma_{i_1}Q$.
Then from ~\cite{BGP},
\[Q=\sigma_{i_n}\cdots \sigma_{i_1} Q.\] 
Define the Coxeter functor $\Phi^+: \mrm{Rep} (Q) \to
\mrm{Rep} (Q)$ to be 
\[\Phi^+ = \Phi_{i_n}^+\circ \cdots \circ \Phi_{i_1}^+.\] 
 
Similarly, if $i$ is a source in $Q$.
For any $(V, x) \in \mrm{Rep}(Q)$, 
let $\Phi_{i}^-(V, x)=(W, y)$ be a representation of 
$\sigma_{i}Q$, 
where $W_j=V_j$ if $j\neq i$ and 
$W_i$ equals the cokernel of the linear map 
$\sum_{\omega\in \Omega: t(\omega)=i}x_{\omega}:
V_i \to \oplus_{\omega\in \Omega: t(\omega)=i}V_{h(\omega)}$; 
for any $\omega \in \Omega$, 
$y_{\omega}=x_{\omega}$ if $h'(\omega)\neq i$, 
otherwise if $h'(\omega)=i$, $y_{\omega}$ is the composition of the natural maps 
$W_{t'(\omega)} \to \oplus_{\omega\in \Omega:t(\omega)=i}V_{h(\omega)} \to W_i$.
This extends to a functor
\[\Phi_i^-: \mrm{Rep}(Q) \to \mrm{Rep}(\sigma_iQ).\] 
The Coxeter functor 
$\Phi^- :\mrm{Rep}(Q) \to \mrm{Rep}(Q)$ 
is defined to be 
\[\Phi^-=\Phi_{i_1}^- \circ \cdots \circ \Phi_{i_{n}}^-.\]

\subsection{Classification of indecomposable representations}
Assume that $Q$ has no oriented cycles.
By using the two Coxeter functors $\Phi^+$ and $\Phi^-$, 
the representations $M \in \mrm{Ind}(Q)$ are classified into four
classes: preprojective, inhomogeneous regular, homogeneous regular and
preinjective. More precisely, a representation $M \in \mrm{Ind}(Q)$
is called preprojective if  $(\Phi^+)^r M=0$, for $r \gg 0$; 
inhomogeneous regular if $(\Phi^+)^rM\simeq M$ for some $r\geq2$;
homogeneous regular if $(\Phi^+)^rM\simeq M$ for any positive integer
$r$; and 
preinjective if $(\Phi^-)^r M =0$, for $r \gg  0$.
In general, $M\in \mrm{Rep}(Q)$ is preprojective
(resp. inhomogeneous regular, homogeneous regular,  preinjective) if all indecomposable direct summands of
$M$ are preprojective (resp. inhomogeneous regular, homogeneous
regular, preinjective).

Let $S_{i}$ be the simple representation corresponding to the vertex $i$. 
It's a representation $(V, x)$ such that $V_i=k$, $V_j=0$ 
if $j\neq i$ and all linear maps $x_{\omega}$ are 0. 
Note that given a graph, the definition of the simple representation 
works for any orientation of the graph. 
By abuse of notation, we always denote by $S_{i}$ 
the simple representation corresponding to the vertex $i$ regardless of the orientation. 

Fix a sequence $(i_1,\cdots,i_n)$ such that $i_r$ is a sink of the quiver
$\sigma_{i_{r-1}}\cdots\sigma_{i_1}(Q)$. 
Let 
\[ P(i_r)=\Phi_{i_1}^-\circ\cdots\circ\Phi_{i_{r-1}}^-(S_{i_r}) 
\quad \text{and} \quad
I(i_r)=\Phi_{i_n}^+\circ\cdots\circ\Phi_{i_{r+1}}^+(S_{i_r}).\] 
Then we have
\begin{align*}
&M\in \mrm{Ind}(Q)\; \text{ is projective iff}\;
M\simeq P(i_r). \tag{1}\\
&M\in \mrm{Ind}(Q)\; \text{ is injective iff}\;
M\simeq I(i_r). \tag{1'}\\
&M\in \mrm{Ind}(Q)\;\text{ is preprojective iff} \;
M=(\Phi^-)^r\, P(i). \tag{2}\\ 
&M \in \mrm{Ind}(Q)\; \text{is preinjective iff}\; 
M =(\Phi^+)^r\, I(i)\tag{3} 
\end{align*}
Let $\mrm{Reg}(Q)$ be the full subcategory of $\mrm{Rep}(Q)$
whose objects are inhomogeneous regular and
homogeneous representations. Then we have
\begin{align*}
&\mrm{Reg}(Q)\;\text{is an extension-closed 
          full subcategory of}\; \mrm{Rep}(Q). \tag{4}\\
&\Phi^+\;\text{is an autoequivalence on}\; \mrm{Reg}(Q). \;\Phi^-\;\text{ is its inverse.} 
\tag{5}
\end{align*}

The simple objects in $\mrm{Reg}(Q)$ are called regular simple representations. 
For each regular simple representation $R$, there exists a positive integer $r$ such that 
$(\Phi^+)^r\, R=R$. 
We call the smallest one, $p$, the period of $R$ under $\Phi^+$. 
The set $\{R,\cdots,(\Phi^+)^{p-1} \, R\}$ 
is called the $\Phi^+$-orbit of $R$.
Given a $\Phi^+$-orbit of a regular simple representation, 
the corresponding $tube$, say $T$, 
is a set of isoclasses of all indecomposable regular representations whose 
regular composition factors belong to this orbit. 
Let $\mrm{Rep}(T)$ be the full subcategory of $\mrm{Rep} (Q)$ the 
objects of which are direct sums of indecomposable representations in
$T$.
We have the following facts:
\begin{align*}
&\text{Every regular indecomposable representation belongs to a unique
  tube;} \tag{6}\\
&\text{Every indecomposable object in a tube has the same period under
  $\Phi^+$;}\tag{7}\\
&\text{All but finitely many regular simple objects have period
  one;}\tag{8}
\end{align*}
Let $p$ be  the cardinality of the set of isomorphism classes of  
simple objects in $\mrm{Rep}(T)$ and $p_T$ be  the period of $T$.
Then
\begin{align*}
&p=p_T
\tag{9}
\end{align*}
Given any representation $M\in \mrm{Rep}(T)$ with $T$ of period $p$, 
$M$ is called aperiodic if for any $N\in T$ not all the indecomposable representations 
\[
N, (\Phi^+)\, N, \cdots, (\Phi^+)^{p-1}\, N
\]
are direct summands of $M$. 
Given any $x\in E_{V, \Omega}$ such that the representation $(V,x)$ is
aperiodic, we call the $G_V$-orbit $O_x$ of $x$ aperiodic.

The following lemma will be needed in the proof of Proposition ~\ref{generalcases}.

\begin{lem} \label{vanishing}
Let $M$ and $ N\in \mrm{Ind}(Q)$ be one of the following cases.
\begin{enumerate}
\item $M=(\Phi^+)^m I(i_r)$ and $ N=(\Phi^+)^{n} I(i_{r'})$,
for $m> n$ or  $m=n$ and $r\leq r'$.
\item $M=(\Phi^-)^{m}P(i_r)$ and
$N=(\Phi^-)^{n}P(i_{r'})$,
for $m<n$ or $m=n$ and $r\leq r'$.
\item $M$ and $N$ are both regular, but they are not in the same tube.
\item $M$ is nonpreinjective and $ N$ is preinjective.
\item $M$ is preprojective and $N$ is nonpreprojective. 
\end{enumerate}
Then we have
\begin{align*}
&\mrm{Ext}^1(M, N)=0; \tag{i}\\
&\mrm{Hom}(N, M)=0;\quad \text{ if}\; M\; \text{ and}\; N\;\text{ are not isomorphic}.
\tag{ii}\\
\end{align*}
\end{lem}

The Lemma follows from  ~\cite[Chapter 6, Lemma 1]{Crawley-Boevey}.

\subsection{Tubes in a noncyclic quiver}
\label{noncyclicquivers}
Let $Q=(I,\Omega; s,t)$ be an affine quiver other than $\mrm C_p$. 
Let $R=(\oplus_{i\in I}R_i, r)$ be an inhomogeneous regular simple in $\mrm{Ind}(Q)$.
Let $T$ be a tube in $\mrm{Rep}(Q)$, consisting of all indecomposable 
regular representations such that there exists a filtration of subrepresentations 
whose consecutive quotients 
are isomorphic to regular simple representations in  the $\Phi^+$-orbit of $R$.
Let $\mrm{Rep}(T)$ be the full subcategory of $\mrm{Rep}(Q)$ 
consisting of direct sums of indecomposable representations in $T$. 
The $\Phi^+$-orbit of $R$ forms a complete list of simple objects in the subcategory  
$\mrm{Rep}(T)$. Let $p$ be the minimal $ k\geq 1 $ such that $(\Phi^+)^k(R)\cong R$.  
For convenience, we write $R_z=(\oplus_{i\in I}R_{z,i},r_z)$ for $(\Phi^+)^k(R)$ 
if $ [k]=z$ in $\mathbb{Z}/p\mathbb{Z}$, where $[k]$ is the class of $k$ in 
$\mathbb{Z}/p\mathbb{Z}$. 
Note that $ R_{z,i}$ is a vector space and $r_z=(r_{z, \omega}) $ 
with $ r_{z, \omega}:R_{z,t(\omega)}\rightarrow R_{z,h(\omega)}$ is a linear transformation 
for each $ \omega \in \Omega$. 
These regular simples in $ T$  satisfy the following properties.
\begin{enumerate}
\item[(1)] $\mrm{Hom}_Q\;(R_z,R_z)=K$ and $\mrm{Hom}_Q\;(R_z,R_{z'})=0$ 
      if $z\neq z'$;
\item[(2)] $\mrm{Ext}^1_Q(R_z,R_{z-1})=K$ and $\mrm{Ext}^1_Q(R_z,R_{z'})=0$ 
      if $z\neq z'+1$;
\item[(3)] All higher extension groups of the simple objects in $\mrm{Rep}(T)$ vanish.
\end{enumerate}
For each $z\in \mathbb{Z}/p\mathbb{Z}$, 
fix an extension $E_z=(\oplus_{i \in I} E_{z,i}, e_z)$ of $R_z$ by $R_{z-1}$ 
such that the exact sequence
$0\to R_{z-1}\to E_z \to R_z \to 0$ is nonsplit. Notice that 
$E_{z,i}=R_{z-1,i}\oplus R_{z,i}$ for any $i\in I$. 
(In another word, we fix
a basis of the extension group $\mrm{Ext}^1_Q(R_z,R_{z-1})$.)
For each arrow $\omega\in \Omega$, 
the linear map $e_{z,\omega}:E_{z,t(\omega)}\to E_{z, h(\omega)}$
induces a linear map 
\[
l_{z,\omega}:R_{z,t(\omega)}\to R_{z-1,h(\omega)}.
\tag{4}
\]

\subsection{Cyclic quivers}
\label{cyclic}
Let $\mrm C_p$ be a cyclic quiver of $p$ vertices. More precisely, 
$\mrm C_p$ is defined to be the quadruple $\mrm C_p=(I_p, \Omega_p; s_p,t_p)$, where 
\begin{enumerate}
\item[(i)] $I_p=\mathbb{Z}/p\mathbb{Z}$, 
\item[(ii)] $\Omega_p=\{\omega_z|\;z\in \mathbb{Z}/p\mathbb{Z}\}$, 
\item[(iii)] $t_p(\omega_z)=z$ and $h_p(\omega_z)=z-1$ for any $z\in \mathbb{Z}/p\mathbb{Z}$. 
\end{enumerate}
We sometimes write $z\to z-1$ for the arrow $\omega_z$.
For each vertex $z\in \mathbb{Z}/p\mathbb{Z}$, denote by $s_z$ the corresponding 
simple representation, i.e., $s_z$ is the representation whose associated vector
space is $K$ at vertex $z$ and $0$ elsewhere and whose associated linear maps are zero.
For each $\lambda \in K^*$, define a representation $t_{\lambda}$ 
by associating the one dimensional vector space $K$ 
to each vertex and  the identity map to each arrow except the arrow $0 \to p-1$ 
which is associated with the scalar map $\lambda \mrm{Id}$, 
where $\mrm{Id}$ is the identity map. 
The union 
$\{s_z| z\in \mathbb{Z}/p\mathbb{Z}\} \cup \{t_{\lambda}|\lambda \in K^*\}$ 
forms a complete list of pairwise nonisomorphic simple representations 
in $\mrm{Rep}(\mrm C_p)$. 
In particular, the subset $\{s_z|\;z\in \mathbb{Z}/p\mathbb{Z}\}$ 
is a complete list of pairwise nonisomorphic simple objects
in $\mrm{Nil}(\mrm C_p)$. 
Moreover, they satisfy the following properties.
\begin{itemize}
\item[(1)] $\mrm{Hom}_{\mrm C_p}\;(s_z,s_z)=K$ and $\mrm{Hom}_{\mrm
           C_p}(s_z,s_{z'})=0$ 
           if $z\neq z'$;
\item[(2)] $\mrm{Ext}^1_{\mrm C_p}(s_z,s_{z-1})=K$ and 
           $\mrm{Ext}^1_{\mrm C_p}(s_z,s_{z'})=0$ if $z'\neq z-1$;
\item[(3)] $\mrm{Hom}_{\mrm C_p}(t_{\lambda},t_{\lambda})=K$ and 
           $\mrm{Hom}_{\mrm C_p}(t_{\lambda},t_{\lambda'})=0$ if $\lambda \neq \lambda'$;
\item[(4)] $\mrm{Ext}^1_{\mrm C_p}(t_{\lambda}, t_{\lambda})=K$ and 
      $\mrm{Ext}^1_{\mrm C_p}(t_{\lambda}, t_{\lambda'})=0$ if $\lambda \neq \lambda'$;
\item[(5)] $\mrm{Hom}_{\mrm C_p}(s_z, t_{\lambda})=0$ and
      $\mrm{Hom}_{\mrm C_p}(t_{\lambda}, s_z)=0$
      for any $z\in \mathbb Z/p\mathbb{Z}$ and $\lambda \in K^*$;           
\item[(6)] $\mrm{Ext}^1_{\mrm C_p}(s_z,t_{\lambda})=0$ and 
      $\mrm{Ext}^1_{\mrm C_p}(t_{\lambda}, s_z)=0$
      for any $z\in \mathbb Z/p\mathbb{Z}$ and $\lambda \in K^*$;       
\item[(7)] All higher extension groups of the simple objects vanish.
\end{itemize} 

Let $s_{z,l}$ be the indecomposable representation of $\mrm C_p$ such
that its socle is $s_i$ and  its length is $l$ for $z\in I_p$ and
$l\in \mathbb N$. Given any representation $M$ in $\mrm{Rep}(\mrm C_p)$,
$M$ is called aperiodic if for each $l\in \mathbb N$ not all of the
indecomposable representations 
\[
s_{0, l}, s_{1, l}, \cdots, s_{p-1,l}
\]
are direct summands of $M$.

Given any $x \in E_{\mathbb V, c_p}$, 
if the representation $(V,x)$ is aperiodic, then we call the
$G_{\mathbb V}$-orbit $O_x$ of $x$ aperiodic.

\section{Hall functors and Hall morphisms}

\subsection{The Hall functor $F$}
Let $T$ be a tube (Section \ref{noncyclicquivers}) in $\mrm{Rep}(Q)$. Recall
that $R_z=(\oplus_{i\in I}R_{z,i},r_z)$ ($z\in \mathbb Z/p\mathbb Z$)
are the pairwise nonisomorphic simple objects in $\mrm{Rep}(T)$.
$l_{z,\omega}:R_{z,t(\omega)}\to R_{z-1,h(\omega)}$
($\omega\in\Omega$) is a linear map defined Section
~\ref{noncyclicquivers}.

For any representation $(\mathbb V,\theta)\in \mrm{Rep}(\mrm C_p)$
(Section ~\ref{cyclic}), 
define a representation 
\[F(\mathbb V,\theta)=(F(\mathbb V),F(\theta))\in \mrm{Rep}(Q)\] 
by
\begin{enumerate}
\item[(i)] $F(\mathbb V)_i=\oplus_{z\in \mathbb{Z}/p\mathbb{Z}}\mathbb V_z\otimes R_{z,i}$ 
          for any $i\in I$, 
\item[(ii)] $F(\theta)_{\omega}=\sum_{z\in \mathbb{Z}/p\mathbb{Z}}
          (\I_z \otimes r_{z,\omega}+\theta_{z\to z-1}\otimes l_{z,\omega})$,
          for each arrow $\omega\in \Omega$.
\end{enumerate}
(Here $\I_z: \mathbb{V}_z\rightarrow \mathbb{V}_z$ in (ii) is the identity map.)
Note that $F(\theta)_{\omega}$ is a linear map for $F(\mathbb
V)_{t(\omega)}$ to $F(\mathbb V)_{h(\omega)}$.
This map extends to 
a functor $F:\mrm{Rep}(\mrm C_p)\to \mrm{Rep}(Q)$.
By construction, we have 
\begin{lem}
 \begin{itemize}
  \item [(a)]  $F$ is an exact functor.
  \item [(b)]  $F(s_z)=R_z$ for all $z \in \mathbb{Z}/p\mathbb{Z}$.
  \item [(c)]  $F(\mathbb V,\theta)\in \mrm{Rep}(T)$ 
              for any $(\mathbb V,\theta) \in \mrm{Nil}(\mrm C_p)$.
  \item [(d)]  $F(t_{\lambda})$ is a homogeneous regular simple for any $\lambda \in K^*$.
 \end{itemize}
\end{lem}

\begin{proof}
From the construction, $F$ is exact on each vector space. To show (a), one only
needs to check that $F$ is a functor which is straightforward
verified. (b) can be checked directly.  

For any $(\mathbb V,\theta)\in \mrm{Nil}(\mrm C_p)$, it is well-known that 
there exists a sequence of subrepresentations
\[(\mathbb V,\theta)\supseteq (\mathbb V^1,\theta^1)\supseteq \cdots 
\supseteq (\mathbb \mathbb V^n,\theta^n)=0,\]
such that $(\mathbb V^l,\theta^l)/(\mathbb V^{l+1},\theta^{l+1})$ are simple representations
of type $s_1, \cdots s_p$. 
From this and the construction of $F$, 
there is a sequence of subrepresentations of $F(\mathbb V,\theta)$:
\[F(\mathbb V,\theta)\supseteq F(\mathbb V^1,\theta^1)\supseteq \cdots 
\supseteq F(\mathbb V^n,\theta^n)=0,\] 
satisfying $F(\mathbb V^l,\theta^l)/F(\mathbb
V^{l+1},\theta^{l+1})\simeq R_z$ by (a)  
for some $z\in \mathbb{Z}/p\mathbb{Z}$. Therefore, (c)  follows from (b).

To prove (d), 
when $Q$ is of type $\tilde{A}_n$, 
the functor $F$ coincides with the functor $G$ in [Li].
From the construction of $G$, $F(t_{\lambda})$ is a homogeneous regular simple.
Furthermore, $\{F(t_{\lambda})|\lambda \in K^*\}$ is a complete list of pairwise nonisomorphic
homogeneous regular simples in $\mrm{Rep}(Q)$.
When $Q$ is of type $\tilde{D}_n$ or $\tilde E_m$ ($m=6, 7$ and $8$), 
we choose the orientation $\Omega'$ given in [DR, P. 46-49]. 
We denote by $Q'$ the corresponding quiver.
By direct computation, $F(t_{\lambda})$ is a homogeneous regular simple. 
Furthermore, $\{F(t_{\lambda})|\lambda \in K^*\}$ is a complete list of pairwise nonisomorphic
homogeneous regular simples in $\mrm{Rep}(Q')$. 
On the other hand,
by [BGP], there exists a sequence $i_1, \cdots, i_k$ of vertices in $I$ such that 
it is $(+)$-accessible with respect to $Q$ and 
$\sigma_{i_k}\sigma_{i_{k-1}}\cdots \sigma_{i_1} Q=Q'$. 
Denote by $\Psi^+:\mrm{Rep}(Q)\to \mrm{Rep}(Q')$ the composition of 
the corresponding reflection functors $\Phi_i^+$. 
We write $R_z'$ for $\Psi^+(R_z)$, for any $z\in \mathbb Z/p\mathbb Z$.
Denote by $T'$ the tube in $\mrm{Rep}(Q')$ generated by $R_z'$, 
for all $z\in \mathbb Z/p\mathbb Z$. 
Define the functor $F': \mrm{Rep}(\mrm C_p)\to \mrm{Rep}(Q')$ as the functor $F$ by
replacing $R_z$ by $R_z'$.
We then have  
\[F'\simeq \Psi^+ \circ F.\]
In fact,
it suffices to prove the case when $\sigma_i Q=Q'$. 
In this case it can be verified directly by the construction of 
the functors $F$, $F'$ and $\Phi_i^+$.

Note that $\Psi^+$ sends homogeneous regular simples to homogeneous regular simples. 
This proves (d).
\end{proof}

We denote by  $\mrm{HT}$  the full subcategory of $\mrm{Rep}(Q)$ generated 
by $R_z$ and $F(t_{\lambda})$, for all $z\in \mathbb Z/p\mathbb Z$, $\lambda\in K^*$. 
From the above Lemma, $F$ induces a functor from $\mrm{Rep}(\mrm C_p)$ to $\mrm{HT}$, 
still denoted by $F$. We have   
\begin{lem} \label{lem:2.6}
The induced functor $F:\mrm{Rep}(\mrm C_p) \to \mrm{HT}$ is a categorical equivalence. 
\end{lem}
\begin{proof}

By  Lemma 2.5 (a), $F$ induces maps:
\begin{enumerate}
\item  $F_z: \mrm{Ext}^1_{\mrm C_p}(s_z,s_{z-1}) \to \mrm{Ext}^1_Q(F(s_z),F(s_{z-1}))$ 
       for any $z\in\mathbb{Z}/p\mathbb{Z}$.
\item  $F_{\lambda}:\mrm{Ext}^1_{\mrm C_p}(t_{\lambda},t_{\lambda})\to 
       \mrm{Ext}^1_Q(F(t_{\lambda}),F(t_{\lambda}))$ for any $\lambda \in K^*$.
\end{enumerate}
hey are all injective and $K$-linear. So $F_z$ and $F_{\lambda}$ are bijective, 
for any $z\in \mathbb{Z}/p\mathbb{Z}$ and $\lambda \in K^*$.
We then have the equivalence by the following Lemma stated in
\cite[P. 129]{GRK}.
\begin{lem}
%$\mbf{Lemma.}$  {\em
Let $E: \mathscr B \to \mathscr C$ be an exact functor between two abelian aggregates 
whose objects have finite Jordan-H\"older series. Then $E$ is an equivalence if and 
only if the following two conditions are satisfied:
\begin{itemize}
\item[(i)] 
$E$ maps simples onto simples and induces a bijection between their sets of isoclasses.

\item[(ii)] For all simples $S, T\in \mathscr B$, the map 
$\mrm{Ext}^i_{\mathscr B}(S,T)\to \mrm{Ext}^i_{\mathscr C}(ES,ET)$ induced by $E$ is bijective
for $i=1$ and injective for $i=2$.
\end{itemize}
\end{lem}
\end{proof}

\noindent
{\bf Remark.}
The functor $F$ is a Hall functor in \cite{FMV}. 
A prototype of the functor $F$ can be found in [GRK].
The restriction $F:\mrm{Nil}(\mrm C_p) \to \mrm{Rep}(T)$ is equivalent.

\subsection{The Hall morphism $F$} 
\label{sec:2.7}
The assignment $(\mathbb V, \theta)\mapsto (F(\mathbb V),F(\theta))$
in the above subsection gives a morphism of varieties:
%\begin{list}%\item[(e).] 
\[ \tag{a} F:  E_{\mathbb V, \Omega_p}\to  E_{F(\mathbb V),\Omega}, \;\;\theta \mapsto F(\theta).\]
Note that both $E_{\mathbb V,\Omega_p}$ and $E_{F(\mathbb V), \Omega}$
are $K$-vector spaces and the map is affine linear ( a linear map plus
a constant map).

We set $V=F(\mathbb V)$. 
By definition, $V$ is an $I$-graded $K$-vector space.
Denote by $E_1$ the image of $E_{\mathbb V,\Omega_p}$ under $F$. 
(Note that $E_1$ is a translate of a vector subspace of $E_{V, \Omega}$.)
Denote by $E_2$ the set of all elements $x \in E_{V,\Omega}$ such that
$x$ is in the same $G_V$-orbit of some element in $E_1$. Clearly,
\[\tag{b}  E_1\subseteq E_2\; \text{ and } \;E_2 \; \text{ is a $G_V$-stable subvariety
of $E_{V, \Omega}$}.\]
Moreover,  the assignment $(\mathbb V, \theta)\mapsto (F(\mathbb V),F(\theta))$
induces an algebraic group homomorphism defined by
\[\tag{c}  F: G_{\mathbb V} \to G_V, \;\; g \mapsto g \otimes 1, \; \text{ 
   for any } \;g \in G_{\mathbb V},\] 
where
$(g \otimes 1)_i=\oplus_{z\in \mathbb{Z}/ p \mathbb{Z}}\, g_z \otimes \I_{z,i} :
\oplus_z \mathbb V_z \otimes R_{z,i} \to \oplus_z \mathbb V_z \otimes R_{z,i}$ 
for all $i\in I$, with $\I_{z,i}: R_{z,i}\rightarrow R_{z,i}$ being the identity map.
For simplicity, we write $H_V$ for $F(G_{\mathbb V})$. 
It is an algebraic subgroup of $G_V$.
Note that $G_{\mathbb V} \simeq H_V$. 
Also there is an action of $H_V$ on $E_1$ induced by the action of
$G_{\mathbb V}$ on $E_{\mathbb V,\Omega_p}$.
By definitions, the action of $G_{\mathbb V}$ on $E_{\mathbb V,\Omega_p}$ is compatible with 
the action of $H_V$ on $ E_1$, i.e.,
\[\tag{d}  
F(g\theta)=F(g)F(\theta),\; \text{ for any }\; g\in G_{\mathbb V} \;
\text{ and }\theta\in E_{\mathbb V,\Omega_p}.\]
Furthermore, we have
\[ \tag{e} \;  
\text{The map $F: E_{\mathbb V, \Omega_p} \to E_1$ is an isomorphism of varieties.}\]
\[\tag{f} \text{ For  $x, x_1\in E_1$,  $\{\xi \in G_V\;|\; \xi
x=x_1\} \subseteq  H_V$. In particular,  $ H_V=\mrm{Stab}_{G_V}(F(0))$.}\]
In fact, 
$x=F(\theta)$ and $x_1=F(\theta_1)$, 
for some $\theta, \theta_1\in E_{\mathbb V,\Omega_p}$.
So $\xi F(\theta)=F(\theta_1)$.
That is
$\xi \in \mrm{Hom}_Q((V,F(\theta)), (V,F(\theta_1)))$. 
Since the functor $F:\mrm{Rep}(\mrm C_p)\to \mrm{HT}$ is a
categorical equivalence, 
it is then fully faithful. So 
$\mrm{Hom}_Q((V,F(\theta)), (V,F(\theta_1)))=
F(\mrm{Hom_{C_p}}((\mathbb V,\theta), (\mathbb V,\theta_1)))$. 
Therefore, $\xi=F(g)$ for some $g\in G_{\mathbb V}$.

Let 
\[
\mrm{GL}(V_i)^0=\prod_{z\in \mathbb{Z/Z}p} \mrm{GL}(\mathbb V_z\otimes R_{z,i})
\quad 
\text{and}
\quad
G_V^0=\prod_{i\in I} \mrm{GL}(V_i)^0.
\]
Then $G_V^0$ is a Levi subgroup of the reductive group $G_V$.
Let 
\begin{equation*}
\begin{split}
E_{V,\Omega}^0=\oplus_{\omega\in \Omega} \oplus_{z\in \mathbb{Z/Z}p} 
\mrm{Hom}_K(\mathbb V_z\otimes R_{z,t(\omega)}, \mathbb V_z\otimes R_{z, h(\omega)}),\\
E_{V,\Omega}^1=\oplus_{\omega\in \Omega} \oplus_{z\in \mathbb{Z/Z}p} 
\mrm{Hom}_K(\mathbb V_z\otimes R_{z,t(\omega)}, \mathbb V_{z-1}\otimes R_{z-1, h(\omega)}).
\end{split}
\end{equation*}
By definition
\[
G_V^0E_1 \subseteq E_{V,\Omega}^0\oplus E_{V,\Omega}^1.
\]
The map $\phi: G_V^0 E_1 \to G_V^0F(0) \quad (gF(\theta)\mapsto gF(0))$ is then a 
restriction of the second projection $E_{V,\Omega}^0\oplus E_{V,\Omega}^1 \to E_{V,\Omega}^0$.
Hence $\phi$ is a morphism of varieties.
By ~\cite[Lemma 4 in 3.7]{slodowy} and the fact that $G_V^0/H_V\simeq G_V^0 F(0)$, we have
\[
\label{k}
G_V^0\times ^{H_V}E_1 \simeq G_V^0 E_1.
\tag{g}
\]

Define an action of $H_V$ on $G_V \times  E_1$ by
$(\zeta, x).\xi=(\zeta\xi, \xi^{-1}x)$ for any 
$(\zeta,x)\in G_V \times  E_1$ and $\xi \in H_V$.
Define a map $\tau: G_V \times E_1 \to E_2$ by $(\zeta, x) \mapsto \zeta.x$ 
for any $(\zeta, x)\in G_V \times  E_1 $. 
By (j), the morphism $\tau$ is an $H_V$-orbit map, i.e., all the
fibres of $\tau$ over   $x\in E_2$ are $H_V$-orbit. 
Thus it induces a bijective morphism of varieties 
\[
\bar\tau: G_V\times^{H_V} E_1 \to E_2.
\]
Note that $H_V \subseteq G_V^0$. By ~\cite{slodowy} and (\ref{k}), we have
\[
G_V\times^{H_V}E_1 \simeq G_V\times^{G_V^0}(G_V^0\times ^{H_V} E_1)\simeq 
G_V\times ^{G_V^0} (G_V^0E_1).
\]
Thus $\bar\tau$ induces a bijective morphism of varieties 
\[
\bar\tau: G_V\times^{G_V^0}(G_V^0E_1) \to E_2.
\tag{h}
\]
Let $U^+$ (resp. $U^-$) be the block upper (resp. lower) triangle matrices 
with respect to $G_V^0$ in $G_V$.
$U^+$ and $U^-$ are unipotent radicals of opposite parabolics with
$G_V^0$ as Levi subgroup. Then $U^+\times G^0_V \times U^-$ is an affine open
subvariety in $G_V$. Thus $U:=U^+\times U^- \times G_V^0E_1$ is an affine open subvariety in 
$G_V\times^{G_V^0} G_V^0E_1$. By restricting to the affine open subvariety
$U$, the morphism $\bar\tau$ becomes an isomorphism onto its image by
direct matrix computations (the inverse of $\bar \tau$ is algebraic). Now that 
$\{gU\,|\, g\in G_V\}$ is an open cover of $G_V\times^{G_V^0}(G_V^0E_1)$, we have
\begin{lem}
\label{geometric-F}
$G_V\times^{H_V} E_1 \simeq E_2$ as $G_V$-varieties.
\end{lem}
In particular,
\[\tag{i}  O_x=G_V \times^{H_V} O'_x,  \]
where $x\in  E_1$, $O_x$ is the $G_V$-orbit of  $x$ in
$ E_2$ and $O'_x$ is the $H_V$-orbit of $x$ in $ E_1$.

\section{The canonical basis}
In this section, we recall Lusztig's geometric realization of the
canonical basis of $\mbf U^-$.
\subsection{Notations}
We fix some notations, most of them are consistent with the notations
in ~\cite{lusztig3}.

Fix a prime $l$ that is invertible in $K$.
Given any algebraic variety $X$ over $K$, 
denote by $\mathcal{D}_c^b(X)$ the bounded derived category of complexes 
of $l$-adic sheaves on $X$ (\cite{BBD}).
Let $\mathcal M(X)$ be the full subcategory of $\mathcal D_c^b(X)$
consisting of all perverse sheaves on $X$ (\cite{BBD}).

Let $G$ be a connected algebraic group. Assume that $G$ acts on $X$
algebraically. Denote by $\mathcal D_G(X)$ the full subcategory of
$\mathcal D_c^b(X)$ consisting of all $G$-equivariant complexes over
$X$.
Similarly, denote by $\mathcal M_G(X)$ the full subcategory of
$\mathcal M(X)$ consisting of all $G$-equivariant perverse sheaves (\cite{lusztig3}).

Let $\bar{\mathbb Q}_l$ be an algebraic closure of the field of $l$-adic numbers.
By abuse of notation, denote by  $\bar{\mathbb Q}_l=(\bar{\mathbb Q}_l)_X$ 
the complex concentrated on degree zero, corresponding to
the constant $l$-adic sheaf over $X$. 
For any complex $K \in \mathcal D_c^b(X)$ and $n\in \mathbb Z$, 
let $K[n]$ be the complex such that $K[n]^i=K^{n+i}$ and the
differential is multiplied by a factor $(-1)^n$.
Denote by $\mathcal M(X)[n]$ the full subcategory of 
$\mathcal D_c^b(X)$ objects of which are of the form $K[n]$ with 
$K\in \mathcal M(X)$.
For any $K\in \mathcal D_c^b(X)$ and $L\in \mathcal D_c^b(Y)$, denote
by $K\boxtimes L$ the external tensor product of $K$ and $L$ in
$\mathcal D_c^b(X\times Y)$.

Let $f: X\to Y$ be a morphism of varieties, denote by
$f^*: \mathcal D_c^b(Y) \to \mathcal D_c^b(X)$ and 
$f_!: \mathcal D_c^b(X) \to \mathcal D_c^b(Y)$ 
the inverse image functor and the direct image functor with compact support, respectively.

If $G$ acts on $X$ algebraically and $f$ is a principal $G$-bundle,
then $f^*$ induces a functor (still denote by $f^*$) of equivalence
between $\mathcal M(Y)[\dim G]$ and $\mathcal M_G(X)$.
Its inverse functor is denoted by $f_{\flat}: \mathcal M_G(X) \to
\mathcal M(Y)[\dim G]$ (\cite{lusztig3}).

\subsection{Lusztig's induction functor}
\label{inductionfunctor}
We recall the (geometric) definition of the canonical basis
from ~\cite{lusztig1} (or ~\cite{lusztig3}).
Let $V=\oplus_{i\in I} V_i$ be an $I$-graded $K$-vector space. 
We  have the following data:
\begin{enumerate}
\item The set $\mathcal{X}_{|V|}$. It consists of all sequences 
      $\boldsymbol{\nu}=(\nu^1,\cdots, \nu^n)$ 
      such that $\sum_{m=1}^{n}\nu^m=|V|$ and 
      $\nu^m_{s(\omega)}\cdot \nu^m_{t(\omega)}=0$, 
      for any $\omega \in \Omega$ 
      and $1\leq m \leq n$.
\item The variety $\mathcal{F}_{\boldsymbol{\nu}}$,  for any
      $\boldsymbol{\nu}=(\nu^1,\cdots,\nu^n)\in \mathcal{X}_{|V|}$. 
      It consists of all sequences, 
      $(V=V^0\supseteq V^1\supseteq \cdots \supseteq V^n=0)$,
      of $I$-graded subspaces of $V$ such that $|V^m/V^{m+1}|=\nu^{m+1}$, 
      for $0\leq m \leq n-1$. 
\item The variety $\tilde{\mathcal{F}}_{\boldsymbol{\nu}}$, for any
      $\boldsymbol{\nu}=(\nu^1,\cdots,\nu^n)\in \mathcal{X}_{|V|}$. 
      It consists of all pairs $(x,\mbf f)$, 
      where $x\in E_{V,\Omega}$ and $\mbf f\in \mathcal{F}_{\boldsymbol{\nu}}$, 
      such that $\mbf f$ is $x$-stable.
      (Here $\mbf f$ is $x$-stable means that for any vector subspace,
      $V^m$ in $\mbf f$, 
      $x_{\omega}(V^m_{s(\omega)})\subseteq V^m_{t(\omega)}$, for all $\omega \in \Omega$.)
\end{enumerate}
The first projection 
$\boldsymbol{\pi_{\nu}}: \tilde{\mathcal F}_{\boldsymbol{\nu}} \to E_{V, \Omega}$
then induced a right derived functor
\[(\boldsymbol{\pi_{\nu}})_!:\mathcal{D}_c^b(\tilde{\mathcal{F}}_{\boldsymbol{\nu}})
\to \mathcal{D}_c^b(E_{V,\Omega}).\]
Note that 
$\boldsymbol{\pi_{\nu}}$ is proper and $\tilde{\mathcal F}_{\boldsymbol{\nu}}$ is smooth.
By the Decomposition theorem in ~\cite{BBD},
$(\boldsymbol{\pi_{\nu}})_!(\bar{\mathbb Q}_l)$ is a semisimple complex in 
$\mathcal{D}_c^b(E_{V,\Omega})$. 
Moreover, $(\boldsymbol{\pi_{\nu}})_!(\bar{\mathbb Q}_l)$ is $G_V$-equivariant since 
$\boldsymbol{\pi_{\nu}}$ is $G_V$-equivariant. 

Let $\mathcal{B}_V$ be the set consisting of all isomorphism classes of 
simple $G_V$-equivariant perverse sheaves on $E_{V,\Omega}$ 
that are in the direct summand of the semisimple complex
$(\boldsymbol{\pi_{\nu}})_!(\bar{\mathbb Q}_l)$ 
(up to shift) for some $\boldsymbol{\nu}\in \mathcal{X}_{|V|}$.
By abuse of language, we say ``a complex is in $\mathcal B_V$''
instead of ``the isomorphism classes of a complex is in $\mathcal B_V$''. 

\begin{example}
When $|V|=i$, $E_{V,\Omega}$ consists of a single point. The isomorphism classes of 
the complex $\bar{\mathbb Q}_l$ is the only element in $\mathcal B_V$. We denote it by 
$F_i^{[1]}$. 
\end{example}

Let $\mathcal{Q}_V$ be the full subcategory of $\mathcal D_c^b(E_{V,\Omega})$
consisting of all complexes on $E_{V,\Omega}$  
isomorphic to a direct sum of shifts of finitely many complexes in $\mathcal{B}_V$.

Let $W\subseteq V$ be an $I$-graded $K$-subspace of $V$ and $T=V/W$.
For any $x\in E_{V,\omega}$ such that
$x_{\omega}(W_{t(\omega)})\subseteq W_{h(\omega)}$ for all $\omega \in
\Omega$, we call that  $W$ is $x$-stable and 
it then induces elements $x_W$ and $x_{T}$ in $E_{W,\Omega}$ and
$E_{T,\Omega}$, respectively. 

Define $E''$ to be the variety consisting of all pairs $(x, V')$, where
$x\in E_{V,\Omega}$ and $V'$ is an $I$-graded subspace of $V$, such that
$|V'|=|W|$ and $V'$ is $x$-stable.

Define $E'$ to  be the variety consisting of all quadruples 
$(X, V'; R', R'')$ where $(x, V')$ is in $E''$, and $R': V'\to W$ and
$R'': V/V' \to T$ are $I$-graded linear isomorphisms.
  
Consider the following diagram
\[
\label{diagram}
\begin{CD}
E_{T,\Omega}\times E_{W,\Omega} @<p_1<< E' @>p_2>> E'' @>p_3>> E_{V,\Omega},
\end{CD}
\tag{*}
\]
where the maps are defined as follows.

$p_3: (x, V') \mapsto x$, $p_2:(x, V'; R', R'') \mapsto (x,V')$, and
$p_1: (x,V';R',R'') \mapsto (x', x'')$, where 
$x'_{\omega}=R'_{h(\omega)} x_{V'} (R'_{t(\omega)})^{-1}$ and
$x''_{\omega}=R''_{h(\omega)} x_{V/V'} (R''_{t(\omega)})^{-1}$ 
for all $\omega \in \Omega$.

Note that $p_3$ is proper, $p_2$ is a principal $G_T\times G_W$-bundle
and $p_1$ is smooth with connected fibres.

From (\ref{diagram}), we can form a functor
\[
(p_3)_! (p_2)_{\flat} (p_1)^*: \mathcal D_c^b(E_{T,\Omega}\times
E_{W,\Omega}) \to \mathcal D_c^b(E_{V,\Omega}).
\]
We write $K\star L:=(p_3)_! (p_2)_{\flat} (p_1)^*(K\boxtimes L)$ for any 
$K\in \mathcal D_c^b(E_{T,\Omega})$ and $L\in \mathcal D_c^b(E_{W,\Omega})$.
\begin{lem}
\label{induction}
\begin{enumerate}
\item $K\star L \in \mathcal Q_V$
      for any $K\in \mathcal Q_T$ and $L\in \mathcal Q_W$.
\item $(K\star L) \star M=K \star (L \star M)$ for any $K\in \mathcal
      Q_T$, $L\in \mathcal Q_W$ and $M\in \mathcal Q_U$.
\item $(\boldsymbol{\pi_{\nu}})_!(\bar{\mathbb Q}_l)= 
      (\boldsymbol{\pi}_{(\nu^1)})_! (\bar{\mathbb Q}_l)
      \star \cdots \star 
       (\boldsymbol{\pi}_{(\nu^n)})_! (\bar{\mathbb Q}_l)$
      where $(\nu^m)$ are sequences with only one entry $\nu^m$ for $m=1,\cdots,n$.
\end{enumerate}
\end{lem}

See ~\cite{lusztig1} for a proof.

By Lemma ~\ref{induction} (2), the expression $K\star L \star M$ makes no confusion.

Let $\boldsymbol{\nu}=(\nu_1, \cdots, \nu_n)$ be a sequence of elements in $\mathbb N[I]$.
(($\nu_1, \cdots, \nu_n$) need not satisfy the conditions 
in Section ~\ref{inductionfunctor} (1).)
Let $V$ be an $I$-graded $K$-space of dimension $\sum_{m=1}^n \nu_m$.
Assume that $V^{(m)}$ ($m=1,\cdots, n$) are $I$-graded subspaces of $V$ such that 
$V=\oplus_{m=1}^n V^{(m)}$.

Define $F''$ to be the variety consisting of all pairs $(x, V^{\bullet})$ 
where $x\in E_{V,\Omega}$ and $V^{\bullet}$ is a flag of type $\boldsymbol{\nu}$ such that 
$V^{\bullet}$ is $x$-stable.

Define $F'$ to be the variety consisting of all triples $(x, V^{\bullet}; \mbf g)$
where $(x, V^{\bullet})$ is in $F'$ and $\mbf g$ is a sequence of linear isomorphisms
$(g_m: V^{m-1}/V^m\to V^{(m)}\;|\; m=1, \cdots, n)$. 
(Here $V^{\bullet}=(V=V^0\supseteq V^1 \supseteq \cdots \supseteq V^n=0)$.)

Consider the following diagram
\[
\label{n-sequence}
\begin{CD}
E_{V^{(1)}, \Omega}\times \cdots \times E_{V^{(n)},\Omega} @<q_1<< F' @>q_2>> F'' @>q_3>> E_{V,\Omega}
\end{CD}
\tag{**}
\]
where the maps are defined by $q_3: (x, V^{\bullet}) \mapsto x$, 
$q_2: (x, V^{\bullet},\mbf g)\mapsto (x, V^{\bullet})$, and 
$q_1: (x, V^{\bullet}, \mbf g)\mapsto (x^{(1)}, \cdots, x^{(n)})$ with
$x^{(m)}_{\omega}= (g_m)_{h(\omega)} x_{V^{m-1}/V^m} (g_m)_{t(\omega)}^{-1}$ for $m=1\cdots,n$.

Similar to $p_1, p_2$ and $p_3$, 
the morphisms $q_1, q_2$ and $q_3$ are smooth with connected fibres, principal
$G_{V^{(1)}}\times \cdots \times G_{V^{(n)}}$-bundle and proper, respectively.
Note that when $n=2$, diagram (\ref{n-sequence}) coincides with diagram (\ref{diagram}).

\begin{lem}
\label{n-mult}
$K^{(1)}\star \cdots \star K^{(n)}
=(q_3)_! (q_2)_{\flat} (q_1)^* (K^{(1)}\boxtimes \cdots \boxtimes K^{(n)})$ 
for any $K^{(m)}\in \mathcal Q_{V^{(n)}}$ where  $m=1, \cdots,n$.
\end{lem}
\begin{proof}
The statement follows from Lemma ~\ref{induction} when $n=2$. 
When $n> 2$, the statement can be proved by  induction.
\end{proof}

\subsection{Lusztig's algebras}
\label{algebra}
Let $\mathcal{K}_V=\mathcal K(\mathcal Q_V)$ be the Grothendieck group of the category
$\mathcal Q_V$, i.e., it is   
the abelian group with one generator $\langle L\rangle$ for each isomorphism class
of objects in $\mathcal{Q}_V$ with relations: 
$\langle L\rangle+\langle L'\rangle=\langle L''\rangle$ 
if $ L'' $ is isomorphic to $L\oplus L'$.

Let $v$ be an indeterminate. Set $\mathbb A=\mathbb Z[v,v^{-1}]$.
Define an $\mathbb A$-module structure on $\mathcal K_V$ by 
$v^n\langle L\rangle =\langle L[n]\rangle $ 
for any generator $\langle L\rangle\in\mathcal{Q}_V$ and 
$n\in \mathbb Z$.
From the construction, 
it is a free $\mathbb A$-module with basis $\langle L\rangle$
where $\langle L\rangle$ runs over $\mathcal{B}_V$.

From the construction, we have 
$\mathcal{K}_V \cong \mathcal{K}_{V'}$,
for any $V$ and $V'$ such that $|V|=|V'|$. 
For each $\nu\in \mathbb N[I]$, fix an $I$-graded vector space $V$ of dimension $\nu$.
Let 
\[
\mathcal{K}_{\nu}=\mathcal{K}_V,\quad 
\mathcal{K}=\oplus_{\nu\in \mathbb{N}[I]}\mathcal{K}_{\nu}\quad
\text{and} \quad
\mathcal{K}_Q=\mathbb{Q}(v)\otimes_{\mathbb A]}\mathcal{K}.
\]
Also let
\[
\mathcal B_{\nu}=\mathcal B_V
\quad \text{and} \quad
\mathcal B_Q=\cup_{\nu\in \mathbb N[I]} \mathcal B_{\nu}.
\]
For any $\alpha, \beta\in \mathbb N[I]$, the operation $\star$ induces an 
$\mathbb A$-linear map
\[
\star: \mathcal K_{\alpha} \otimes_{\mathbb A} \mathcal K_{\beta} 
\to \mathcal K_{\alpha +\beta}.
\]
By adding up these linear maps, we have a linear map
\[
\star: \mathcal K \otimes_{\mathbb A} \mathcal K \to \mathcal K.
\] 
Similarly, the operation $\star$ induces a $\mathbb Q(v)$-linear map
\[
\star:\mathcal K_Q \otimes_{\mathbb Q(v)} \mathcal K_Q \to \mathcal K_Q.
\]
\begin{prop}
\begin{enumerate}
\item $(\mathcal K, \star)$ (resp. $(\mathcal K_Q, \star)$) is an associative algebra 
      over $\mathbb A$ (resp. $\mathbb Q(v)$).
\item $\mathcal B_Q$ is an $\mathbb A$-basis of $(\mathcal K, \star)$ 
      and a $\mathbb Q(v)$-basis of $(\mathcal K_Q, \star)$. 
\end{enumerate}
\end{prop}
\begin{proof}
The associativity of $\star$ follows from Lemma ~\ref{induction} (2). 
\end{proof}

Define a new $\mathbb A$-linear map 
$\circ: \mathcal K_{\alpha} \otimes \mathcal K_{\beta} \to \mathcal K_{\alpha+\beta}$ 
by 
\[x\circ y=v^{m(\alpha,\beta)} x\star y\] 
where
\[m(\alpha, \beta)=\sum_{i\in I}\alpha_i \beta_i+
\sum_{\omega\in \omega}\alpha_{t(\omega)}\beta_{h(\omega)}.\]
This induces a bilinear map 
\[
\circ: \mathcal K\otimes \mathcal K \to \mathcal K.
\]
Then 
\begin{cor}
$(\mathcal K, \circ)$ is an associative algebra over $\mathbb A$.
\end{cor}
The linear map $\circ$ satisfies the associativity due to the fact that $m(-,-)$ is a cocycle.
Similarly, we have an associative algebra $(\mathcal K_Q, \circ)$ over $\mathbb Q(v)$.

\subsection{The canonical basis $\mbf B$ of  the algebra $\mbf U^-$}
\label{canonical}
Given any quiver $Q$, let $c_{ii}=2$ and  
\[c_{ij}=-\#\{ \omega\in \Omega\;|\; 
\{t(\omega), h(\omega)\}=\{i,j\}\} \quad \text{for} \quad i\neq j.\]
$C=(c_{ij})_{i,j\in I}$ is then a symmetric generalized Cartan matrix.
Note that the Cartan matrix $C$ is independent of changes of the orientation of $Q$.

For any $m\leq n\in \mathbb N$, let 
\[
[n]=\frac{v^n-v^{-n}}{v-v^{-1}}, \quad 
[n]^!=\prod_{m=1}^n [m]
\quad
\text{and}
\quad
\bin{n}{m}=\frac{[n]^!}{[m]^![n-m]^!}.
\]

Denote by  $\mbf{U}^-$  the negative part of the quantized
enveloping algebra attached to the Cartan matrix $\mrm{C}$. 
$\mbf U^-$ is the quotient of the free algebra with generators $F_i$, $i\in I$ by
the two-sided ideal generated by 
\begin{equation}
\sum_{p=0}^{1-c_{ij}} (-1)^p \bin{1-c_{ij}}{p} F_i^p F_j F_i^{1-c_{ij}-p}, 
\end{equation}
for $i\neq j\in I$.

Let $F^{(n)}_i=\frac{F_i^n}{[n]^!}$ for all $i\in I$ and $n\in \mathbb N$.
Let $_{\mathbb A}\mbf U^-$ be the $\mathbb A$-subalgebra of $\mbf U^-$
generated by $F^{(n)}_i$ for $i\in I$ and $n\in \mathbb N$.
We have 

\begin{thm}(\cite{lusztig1}, \cite{lusztig3})
\label{iso}
The map $F_i^{(1)} \mapsto F_i^{[1]}$ induces an $\mathbb A$-algebra isomorphism
\[_{\mathbb A}\phi: \,_{\mathbb A} \mbf U^- \to (\mathcal K, \circ)\]
and a $\mathbb Q(v)$-algebra isomorphism
\[
\phi: \mbf U^- \to(\mathcal K_Q, \circ).
\]
\end{thm}

\noindent
{\bf Remark.} See ~\cite{lusztig3} for a more general treatment that
works for any symmetrisable generalized Cartan matrix $C$.

Given another quiver $Q'=(I,\Omega,t',h')$ such that 
\[\{t(\omega), h(\omega)\}=\{t'(\omega),h'(\omega)\}\] for all $\omega\in \Omega$.
From Theorem ~\ref{iso}, the map $F_i^{(1)}\mapsto F_i^{[1]}$ induces 
an $\mathbb Q(v)$-algebra isomorphism
\[
\phi': \mbf U^- \to (\mathcal K_{Q'}, \circ).
\]
We have 
\begin{thm}(\cite{lusztig1}, \cite{lusztig3})
\label{coincide}
$\phi^{-1}(\mathcal B_Q)=(\phi')^{-1} (\mathcal B_{Q'}).$
\end{thm}

\begin{Def}
$\mbf B=\phi^{-1}(\mathcal B_Q)$ is called  the Canonical Basis of $\mbf U^-$.
\end{Def}

For each $Q$, $\mathcal B_Q$ gives a presentation of $\mbf B$.
The main goal of this paper is to describe the elements in $\mathcal B_Q$ by specifying 
their supports and the corresponding local systems when $Q$ is affine.

\section{The description of the elements in $\mathcal B_Q$ via quiver representations}

\subsection{Simple equivariant perverse sheaves}
Note that elements in $\mathcal B_Q$ are isomorphism classes of simple equivariant perverse
sheaves. We give a brief description of simple equivariant perverse sheaves.

Let $X$  be an algebraic variety over $K$ with a connected algebraic group $G$ acting on it.
Let $Y$ be a smooth, locally closed, irreducible $G$-invariant subvariety of $X$ and 
$\mathcal L$ an irreducible, $G$-equivariant, local system on $Y$.
Denote by $j: Y\to X$ the natural embedding. 
\begin{thm} (\cite{BBD}, ~\cite{BL})
The complex 
\[\mrm{IC}(Y, \mathcal L):=j_{!\star}(\mathcal L)[\dim Y]\] 
is a simple $G$-equivariant perverse sheaf on $X$.
Moreover, all simple $G$-equivariant perverse sheaves on $X$ are of this form.
\end{thm}

\subsection{Cyclic quivers}
When the quiver $Q$ is the cyclic quiver $\mrm C_p$ for some $p\in \mathbb N$. 
The description of the elements in $\mathcal B$ is given as follows. 

Let $\mathbb V$ be a  $\mathbb Z/p\mathbb Z$-graded $K$-vector space.
Recall that  a $G_{\mathbb V}$-orbit $O$ in $E_{\mathbb V, \Omega_p}$
is aperiodic if for any $x \in O$, the representation $(\mathbb V, x)$ is aperiodic
(see Section ~\ref{cyclic}). 
Let $\mathcal O_{\mathbb V}^a$ be the set of all aperiodic $G_{\mathbb V}$-orbits 
in $ E_{\mathbb V,\Omega_p}$.
Given $O \in \mathcal O_{\mathbb V}^a$,
let $\mrm{IC}(O, \bar{\mathbb Q}_l)$ be the intersection cohomology complex on 
$E_{\mathbb V,\Omega}$ determined by the subvariety $O$ and the constant sheaf 
$\bar{\mathbb Q}_l$ on $O$.
The assignment $O\mapsto \mrm{IC}(O,\bar{\mathbb Q}_l)$
defines a map $\mathcal O_{\mathbb V}^a\to \mathcal B_{\mathbb V}$. 
Furthermore, we have 

\begin{thm} 
(\cite[5.9]{lusztig2})
\label{cycliccase}
The map  $\mathcal O_{\mathbb V}^a\to \mathcal B_{\mathbb V}$ is bijective.
\end{thm}

\subsection{Noncyclic quivers}
\label{noncyclic}
From now on, we assume that the affine quiver $Q$ is not $\mrm C_p$, 
for any $p\in \mathbb N$. We follow Lusztig's argument in
~\cite[Section 6]{lusztig2}.
We study three special cases in this section. 

First,
given $M \in \mrm{Ind}(Q)$ of dimension vector $\nu$. 
Assume that $M$ is either preprojective or preinjective.
Let $V$ be a $K$-vector space such that $|V|=\nu$. 
Let $O_M$ be the $G_V$-orbit in $E_{V,\Omega}$ corresponding to $M$.
Then we have

\begin{lem} (\cite[Lemma 6.8]{lusztig2}) 
\label{Lemma-1}
$\mrm{IC}(O_M,\bar{\mathbb Q}_l)\in \mathcal B_V$.
\end{lem}

\begin{proof}
Since $M$ is either preprojective or preinjective, its self-extension group
$\mrm{Ext}^1_Q(M,M)=0$, so the corresponding $G_V$-orbit
$O_M$ is open in $E_{V,\Omega}$ (see ~\cite{Crawley-Boevey}). 
Since $E_{V,\Omega}$ is smooth, $\mrm{IC}(O_M,\bar{\mathbb Q}_l)$ 
is the constant sheaf $\bar{\mathbb Q}_l$ on $E_{V,\Omega}$
up to shift (see ~\cite[Lemma 4.3.2]{BBD}). 
Now that $Q$ has no oriented cycles, 
we can order the vertices in $I$  
$i_1,\cdots,i_n$ ($n=|I|$)  in a way 
such that $i_r$ is a source of the full subquiver $Q_r$ 
with vertex set $I-\{i_1,\cdots,i_{r-1}\}$. 
For any $V$ of dimension vector $\nu$, let 
$\boldsymbol{\nu}=(\nu_{i_1}\;i_1,
\cdots,\nu_{i_n}\;i_n)$.
By definition, 
$\mathcal{F}_{\boldsymbol{\nu}}$ consists of  a single flag. 
Also any $x$ in $E_{V,\Omega}$ stabilizes this flag. 
So the first projection    
$\boldsymbol{\pi_{\nu}}:\tilde{\mathcal{F}}_{\boldsymbol{\nu}} \to E_{V,\Omega}$ 
is an isomorphism. 
Therefore, we have 
$(\boldsymbol{\pi_{\nu}})_!(\bar{\mathbb{Q}}_l)[d]=
\bar{\mathbb{Q}}_l[d]\in \mathcal B_V $.
\end{proof}

Second, we assume that $V$ is an $I$-graded $K$-vector space of dimension vector 
$q \delta$, where $\delta\in \mathbb N[I]$ is the minimal positive imaginary root 
of the symmetric Euler form (see Section ~\ref{preliminary}) associated to Q. 
We define two varieties as follows.
\begin{enumerate}
 \item The variety $X(0)$. 
       It is the subvariety of
       $E_{V,\Omega}$ consisting of all elements $x$ such that 
       $(V,x) \simeq R_1 \oplus \cdots \oplus R_q$,
       where $R_1,\cdots,R_q$ are pairwise nonisomorphic homogeneous regular simples.  
 \item The variety $\tilde{X}(0)$.  
       This variety consists of all pairs $\{x, (R_1,\cdots,R_q)\}$,      
       where $x\in X(0)$ and $(R_1,\cdots,R_q)$ is a sequence of  
       representations in $\mrm{Ind}(Q)$,
       such that $(V,x)\simeq R_1\oplus\cdots\oplus R_q$.
\end{enumerate}
Note that the dimension vectors of $R_m$ in (1) have to be $\delta$, 
for all $m\in \{1, \cdots, q\}$ 
and once $x$ is fixed, the set of the representations $R_1,\cdots, R_q$ in (2) 
is completely determined.
Note also that the closure of $X(0)$ equals $E_{V,\Omega}$. 

In fact, by \cite{Ringel2}, we have $\mrm{dim}\; X(0)=\mrm{dim}\; O_x+q$. 
Since 
\[
\mrm{dim} \; G_V-\mrm{dim} \; E_{V,\Omega}
=\, <q\delta, q\delta>\,= 
\mrm{dim}\;\mrm{Hom}_Q((V,x),(V,x))-\mrm{dim}\;\mrm{Ext}^1((V,x),(V,x)), 
\]
we have
\[
\mrm{dim}\; E_{V,\Omega}=\mrm{dim} \; G_V-\mrm{dim}\;\mrm{Hom}((V,x),(V,x))+q.
\]
So $\mrm{dim}\; E_{V,\Omega}=\mrm{dim}\; O_x+q=\mrm{dim} \; X(0)$. 
Therefore $E_{V, \Omega}$ is the closure of $X(0)$.

The first projection 
\[
\pi_1:\tilde{X}(0)\to X(0)
\quad
(x, (R_1,\cdots, R_q)) \mapsto x
\] 
is an $S_q$-principal covering
where $S_q$ is the symmetric group  of $q$ letters. 
$S_q$ acts naturally on $(\pi_1)_{\star}(\bar{\mathbb{Q}}_l)$. 
Given any irreducible representation $\chi$ of $S_q$, 
denote by $\mathscr L_{\chi}$  the irreducible local system corresponding to 
the representation $\chi$ via the monodromy functor (see ~\cite{Iversen}).
Note that $\mathscr L$ is a direct summand of 
$(\pi_1)_{\star}(\bar{\mathbb{Q}}_l)$.

Let $\mrm{IC}(X(0),\mathscr L_{\chi})$ be  the intersection complex on 
$E_{V,\Omega}$ determined by $X(0)$ and $\mathscr L_{\chi}$. 
We then have:

\begin{lem} (\cite[6.10 (a)]{lusztig2})
\label{Lemma-2}
$\mrm{IC}(X(0),\mathscr L_{\chi})\in \mathcal B_V.$
\end{lem}

\begin{proof}
Let $\boldsymbol{\delta}=(\delta,\cdots,\delta)$
such that 
$|\boldsymbol{\delta}|=q\;\delta$. 
Let 
$\boldsymbol{\pi_{\delta}}:\tilde{\mathcal{F}}_{\boldsymbol{\delta}} 
\rightarrow E_{V,\Omega}$ 
be the first projection defined as the morphism $\boldsymbol{\pi_{\nu}}$ 
in Section ~\ref{inductionfunctor}. 
Given any $x\in X(0)$, 
let 
\[\mbf f=(V=V^0\supseteq V^1\supseteq\cdots\supseteq V^n=0)\] 
be a flag in
$\mathcal F_{\boldsymbol{\nu}}$ such that $\mbf f$ is $x$-stable. Then
\[
\tag{1} 
\text{The subrepresentation $(V^m,x)$ is regular, for any $m\in \{0,\cdots,n-1\}$}.
\]
This is because if $x\in X(0)$, $(V,x)$ is a regular representation. 
So $(V^r,x)$ can not have preinjective subrepresentations. 
Now that the dimension vector of $(V^r,x)$ is $\delta$, 
the defect of $(V^r,x)$ is zero. 
Thus, $(V^r,x)$ can not have preprojective subrepresentations.
Therefore, $(V^r,x)$ is regular.

Fix an element   $x\in X(0)$,  
we decompose $V=\oplus_{r=1}^q V(r)$ such that $V(r)$ is $x$-stable and $|V(r)|=\delta$.
By the definition of $X(0)$ and (1), we have $(V,x) \simeq \oplus_r (V(r),x)$.
\[
\tag{2}
\text{ Moreover, this decomposition is unique up to order.} 
\]
In fact, if $V=\oplus_r W(r)$ is another decomposition, 
we can reorder the $W(r)$'s such that the subrepresentations  
$(W(r),x)$ and $(V(r),x)$ are isomorphic, for any $r$. 
Fix an isomorphism $f_r:(W(r),x)\to (V(r),x)$ for each $r$, then
they induce an isomorphism 
\[f:=\sum_r f_r:(V,x)\to (V,x)\] 
satisfying $f(W(r))\subseteq V(r)$, for any $r$.
On the other hand, the composition 
\[
p_{r'} \circ f \circ  i_r: V(r) \overset{i_r}{\to}V \overset{f}{\to} V 
\overset{p_{r'}}{\to} V(r')
\] 
is naturally a homomorphism of representations
in $\mrm{Hom}_Q ((V(r),x),(V(r'),x))$,
where $p_{r'}$ and $i_r$ are natural projection and inclusion, respectively. 
Note that 
\[
\text{$\mrm{Hom}_Q((V(r),x), (V(r'),x))=0$, if $r\neq r'$.}
\]
We have $p_{r'} \circ f \circ  i_r=0$, for $r\neq r'$. 
So, $f(V(r))\subseteq V(r)$. 
But by definition, $f(V(r))\subseteq W(r)$. 
Thus, $V(r)=W(r)$. 
Therefore the decomposition is unique up to order.

From (2), we can define an injective map
$\alpha:\tilde{X}(0)\to \tilde{\mathcal F}_{\boldsymbol{\delta}}$
by $\{x,(R_1,\cdots,R_q)\} \mapsto (x,\mbf f)$,
where $\mbf f$ is the flag 
($V=V^1\supseteq V^2 \supseteq \cdots \supseteq V^{q+1}=0$)
such that $V^r=\oplus_{k=r}^qV(k)$ and $(V(k), x)\simeq R_k$, for $r=1,\cdots,q$.
We then have the following commutative diagram:
\[
\begin{CD}
 \tilde{X}(0) @>\pi_1>>  X(0)\\
 @V\alpha VV             @VVV\\
 \tilde{\mathcal F}_{\boldsymbol{\delta}} @>\boldsymbol{\pi_{\delta}}>>  E_{V,\Omega}\;. 
\end{CD}
\]
Note that  this diagram is Cartesian. So the restriction of
$\boldsymbol{\pi_{\delta}}_!(\bar{\mathbb Q}_l)$ to $X(0)$ is 
$(\pi_1)_!(\bar{\mathbb Q}_l)$. 
Recall that $\mathscr L_{\chi}$ is a direct summand of $(\pi_1)_!(\bar{\mathbb Q}_l)$ 
and $X(0)$ is open in $E_{V,\Omega}$.
So 
\[
\tag{3}
\text{
$\mrm{IC}(X(0),\mathscr L_{\chi})$ is a direct summand of 
$(\boldsymbol{\pi_{\delta}})_!(\bar{\mathbb{Q}}_l)$, up to shift.}
\] 
By Lemma ~\ref{induction} (3),
\[
\tag{4}
(\boldsymbol{\pi_{\delta}})_!(\bar{\mathbb{Q}}_l)
=(\boldsymbol{\pi}_{\delta})_!(\bar{\mathbb Q}_l)\star \cdots 
\star (\boldsymbol{\pi}_{\delta})_!(\bar{\mathbb Q}_l),\]
where $\delta$ is regarded as a sequence with only one entry.
Observe that $(\boldsymbol{\pi}_{\delta})_!(\bar{\mathbb Q}_l)=\bar{\mathbb Q}_l$.
By the proof of Lemma ~\ref{Lemma-1}, they are all in $\mathcal B$ (up to shifts).
From (3) and (4),  $\mrm{IC}(X(0),\mathscr L_{\chi})  \in \mathcal B_V$. 
Lemma ~\ref{Lemma-2} is proved.
\end{proof}

Finally, 
let $T$ be a tube of period $p \neq 1$. 
Let $\mathcal O_{V,T}^a$ be the set of all  aperiodic
$\mrm G_V$-orbits $O_x$ in $E_{V,\Omega}$ (see ~\ref{noncyclicquivers}). 
Given any $O\in \mathcal O_{V,T}^a$, 
denote by  $\mrm{IC}(O,\bar{\mathbb Q}_l)$ the intersection complex on $E_{V,\Omega}$ 
determined by $O$ and the constant local system $\bar{\mathbb{Q}}_l$ on $O$.                 
Then, we have 

\begin{lem}
(\cite[6.9 (a)]{lusztig2})
\label{tube}
$\mrm{IC}(O,\bar{\mathbb Q}_l) \in \mathcal B_V$, for any $O \in \mathcal O_{V,T}^a$.
\end{lem}

\begin{proof}
Fix a regular simple $R$ in $T$, following the construction in Section ~\ref{sec:2.7},
we have an categorical equivalence $F:\mrm{Rep}(\mrm C_p)\to \mrm{HT}$, 
where $\mrm{HT}$ is the full subcategory of $\mrm{Rep}(Q)$ generated by $T$ and 
all the homogeneous regular simples such that $HT$ is closed under
extensions in $\mrm{Rep}(Q)$  and taking kernel and cokernels of
morphisms in $HT$.
Given an element $x$ in $O$, there exists a representation 
$(\mathbb V,\theta)\in \mrm{Rep}(\mrm C_p)$ such that $F(\mathbb V,\theta)\simeq (V,x)$.
In particular, $F(\mathbb V)\simeq V$ as  $K$-vector spaces. 
Consequently, $E_{F(\mathbb V),\Omega}\simeq E_{V,\Omega}$. 
So we can  identify $F(\mathbb V)$ with $V$ and identify 
the $G_{F(\mathbb V)}$-orbit of $F(\theta)$ in $E_{F(\mathbb V),\Omega}$ with $O$
in $E_{V,\Omega}$.

Since $F$ is equivalent and $O$ is aperiodic, 
we have $O_{\theta}$ is aperiodic in $\mrm{Rep}(\mrm C_p)$.
By Theorem ~\ref{cycliccase}, we have 
\[
\tag{1}
\text{
$\mrm{IC}(O_{\theta},\bar{\mathbb Q}_l)\in \mathcal B_{\mathbb V}$.
}
\]
In other words,
\[
\tag{2}
\text{
$\mrm{IC}(O_{\theta},\bar{\mathbb Q}_l)$ is a direct summand of 
$(\pi_{\mbf z})_!(\bar{\mathbb Q}_l)$, (up to shift)
}
\]
for some $\mbf z=(z_1,\cdots,z_n)$, $z_s\in \mathbb Z/p\mathbb Z$.

Define $\mathcal{F}_{\boldsymbol{\nu}}'$ to be the variety consisting of
all flags of the form 
\[F(f)=(F(\mathbb V)\supseteq F(\mathbb V^1)\supseteq
\cdots\supseteq F(\mathbb V^n)),\] 
where
$f=(\mathbb V\supseteq \mathbb V^1\supseteq \cdots \supseteq \mathbb V^n)$ is a flag of type
$\mbf z$.

Define $\tilde{\mathcal F}_{\boldsymbol{\nu}}'$ to be the variety consisting of 
all pair $(x,\mbf f)$, where $x\in E_1$ and  $\mbf f \in \mathcal F_{\boldsymbol{\nu}}'$, 
such that $\mbf f$ is $x$-stable.

Define $\tilde{\mathcal F}''_{\boldsymbol{\nu}}$ to be the variety consisting of 
all pairs $(x,\mbf f)$, where $x\in E_2$ and  $\mbf f\in \mathcal F_{\boldsymbol{\nu}}$ 
such that $\mbf f$ is $x$-stable. Then we have
$\tilde{\mathcal F}''_{\boldsymbol{\nu}}=G_V\times^{H_V}
\tilde{\mathcal F}_{\boldsymbol{\nu}}$

Consider the following commutative diagram
\[
\begin{CD}
\tilde{\mathcal F}_{\mbf z} 
@>\tilde F>> 
\tilde{\mathcal F}'_{\boldsymbol{\nu}}
@>\tilde i_1>>
\tilde{\mathcal F}''_{\boldsymbol{\nu}}
@>\tilde{i}>>
\tilde{\mathcal F}_{\boldsymbol{\nu}}\\
@V\boldsymbol{\pi}_{\mbf z}VV @V\boldsymbol{\pi_{\nu}}'VV 
@V \boldsymbol{\pi_{\nu}}''VV @V \boldsymbol{\pi_{\nu}}VV\\
E_{\mathbb V, \Omega_p} @>F>> E_1 @>i_1>> E_2 @>i>> E_{V,\Omega},
\end{CD}
\]
where the vertical maps are first projections, $i, i_1, \tilde i$ and
$\tilde i_1$ are inclusions, and 
$\tilde F: (\theta, f) \mapsto (F(\theta), F(\mbf f))$.
Note that all squares are Cartesians.

Let $O'$ be the $H_V$-orbit of $F(\theta)$ in $E_1$.
Denote by  $\mrm{IC}(O',\bar{\mathbb Q}_l)$ the intersection complex on
$E_1$ determined by $O'$ and the trivial local system $\bar{\mathbb Q}_l$.

From  Section ~\ref{sec:2.7} (e), 
the statement (2) and the Cartesian square on the left in the diagram above, we
have
\[
\tag{3}
\text{
 $\mrm{IC}(O',\bar{\mathbb Q}_l)$ is a direct summand of 
$(\pi_{\boldsymbol{\nu}}')_!(\bar{\mathbb Q}_l)[d']$, for some $d'$.
}
\]
Let $O''$ be the $G_V$-orbit of $F(\theta)$ in $E_2$. 
(In fact, $O''=O$.)
Denote by $\mrm{IC}(O'', \bar{\mathbb Q}_l)$ be the intersection
complex on $E_2$ determined by $O''$ and $\bar{\mathbb Q}_l$.

By Section ~\ref{sec:2.7} (i) and ~\cite[Theorem 2.6.3]{BL}, the derived functor
\[i_1^*: \mathcal D_{G_V}(E_2) \to \mathcal D_{H_V}(E_1)\]
is a categorical equivalence. In particular,
\[
i_1^*(\mrm{IC}(O'', \bar{\mathbb  Q}_l))=\mrm{IC}(O',\bar{\mathbb Q}_l).
\tag{4}
\]
Since the middle square in the above diagram is Cartesian and
$\boldsymbol{\pi_{\nu}}''$ is proper, by the base change Theorem for
proper morphism (\cite[Theorem 6.2.5]{BBD}), We have 
\[\tag{5} i_1^*(\pi''_{\boldsymbol{\nu}})_!(\bar{\mathbb Q}_l)=
(\pi_{\boldsymbol{\nu}}')_!(\bar{\mathbb Q}_l).\]
Thus by (3), (4) and (5), we have
\[
\tag{6}
\text{
$\mrm{IC}(O'',\bar{\mathbb Q}_l)$ is a direct summand of 
$(\pi''_{\boldsymbol{\nu}})_!(\bar{\mathbb Q}_l)[d'']$, for some $d''$.
}
\]
The right square is Cartesian, so we have 
\[
i^*(\pi_{\boldsymbol{\nu}})_!(\bar{\mathbb Q}_l) 
=(\pi''_{\boldsymbol{\nu}})_!(\bar{\mathbb Q}_l).
\tag{7}\] 
Note that the closure of $E_2$ is
$\boldsymbol{\pi_{\nu}}(\tilde{\mathcal F}_{\boldsymbol{\nu}})$ 
and $E_2$ is open in its closure. By (6) and (7), we have
\[
\text{
$\mrm{IC}(O,\bar{\mathbb Q}_l)$ is a direct summand of 
$(\pi_{\boldsymbol{\nu}})_!(\bar{\mathbb Q}_l)[d]$, for some $d$.
}
\tag{8}
\]
But $(\pi_{\boldsymbol{\nu}})_!(\bar{\mathbb Q}_l)$ is a direct sum of 
simple perverse sheaves from $\mathcal B$ with shifts. Therefore,
$\mrm{IC}(O,\bar{\mathbb Q}_l)\in \mathcal B_V$.  Lemma ~\ref{tube} follows.
\end{proof}

\subsection{General Cases} 
\label{general}
In this section, we study general cases. 

Recall that $\mrm{Ind}(Q)$ is the set of representatives of pairwise
nonisomorphic indecomposable representations of $Q$.

Given any $\nu \in \mathbb N[I]$, denote by $\Delta_{\nu}$  the set
of all pairs $(\sigma,\lambda)$ where
$\sigma: \mrm{Ind}(Q) \to \mathbb N$ is a function and $\lambda$ is
the sequence $(0)$ 
or  a sequence $(\lambda_1,\cdots,\lambda_n)$ of decreasing positive integers
satisfying  the following properties:
\begin{itemize}
  \item [(a)] $\prod_{m=0}^{r-1} \sigma((\Phi^+)^m(R))=0$
                   for any regular representation $R\in
                   \mrm{Ind}(Q)$ of period $r$;
  \item[(b)]  $\sum_{M\in \mrm{Ind}(Q)}\sigma(M)|M| +\sum_{m=1}^n \lambda_m \delta=
                   \nu$ if $\lambda=(\lambda_1,\cdots,\lambda_n)$;
   
  \item[(c)] $\sum_{M\in \mrm{Ind}(Q)}\sigma(M)|M|=\nu$ if $\lambda=(0)$.   
\end{itemize}
From (a), if $R$ is homogeneous, $\sigma(R)=0$. 
From (b) and (c), the function $\sigma$ has finite support.

Given any $(\sigma,\lambda)\in \Delta_{\nu}$, 
fix a $K$-vector space $V$ of dimension vector $\nu$. 
Define the varieties $X(\sigma, \lambda)$ and $\tilde
X(\sigma,\lambda)$, the map $\pi_1$ and 
the irreducible local system $\mathscr L_{\lambda}$ as follows.

If $\lambda=(\lambda_1,\cdots,\lambda_n)$, $X(\sigma,\lambda)$ is the
subvariety of $E_{V,\Omega}$ consisting  of all elements $x$ such that 
\[
(V,x)\simeq \oplus_{M\in \mrm{Ind}(Q)}
M^{\sigma(M)} \oplus R_1 \oplus \cdots  \oplus R_q,\] 
where $ M^{\sigma(M)}$ is the direct sum of $\sigma(M)$ copies of $M$, 
$R_1,\cdots,$ and $R_q$ are pairwise nonisomorphic homogeneous regular simples.         
The variety
$\tilde{X})(\sigma,\lambda)$ is  the variety consisting  of all pairs 
\[
(x, (R_1,\cdots, R_q))\]
where  $x\in X(\sigma, \lambda)$ and
$(R_1,\cdots, R_q)$ is a sequence of homogeneous regular simples in $\mrm{Ind}(Q)$
completely determined by $x$ up to order.
The map $\pi_1: \tilde{X}(\sigma, \lambda) \to X(\sigma, \lambda)$ is the first projection.
Note that the first projection 
$\pi_1:\tilde{X}(\sigma, \lambda) \to  X(\sigma,\lambda)$ is a
$S_q$-principal covering. 
The sequence $\lambda$ determines an irreducible representation
$\chi(\lambda)$
of the symmetric group $S_q$. 
Define $\mathscr L_{\lambda}$ to be 
the direct summand of $(\pi_1)_{\star}(\bar{\mathbb Q}_l)$
corresponding to the irreducible representation of $S_q$ determined by the partition $\lambda$.

If $\lambda=(0)$, the variety $X(\sigma, \lambda)$ is the subvariety
of $E_{V,\Omega}$ consisting of all elements $s$ such that
$
(V,x)\simeq \oplus_{M\in \mrm{Ind}(Q)} M^{\sigma(M)}.
$
$\tilde{X}(\sigma,\lambda)$ is $X(\sigma,\lambda)$.
$\pi_1$ is  the identity map $\tilde{X}(\sigma, \lambda) \to X(\sigma,
\lambda)$.
Denote by $\mathscr L_{\lambda}$ the trivial local system on $X(\sigma,\lambda)$.

Note that when $\lambda=(0)$, the variety $X(\sigma,\lambda)$ is a
$G_V$-orbit in $E_{V,\Omega}$.

Let 
$\mrm{IC}(\sigma,\lambda)=\mrm{IC}(X(\sigma,\lambda),\mathscr L_{\lambda})$
be the simple perverse sheaf on $\mrm E_{V,\Omega}$
determined by $X(\sigma, \lambda)$ and $\mathscr L_{\lambda}$.
\begin{prop}
(\cite[Proposition 6.7]{lusztig2})
\label{generalcases}
$\mrm{IC}(\sigma,\lambda)\in \mathcal B_V$.
\end{prop}

The proof will be given in the next section.

By Proposition ~\ref{generalcases}, the assignment 
$(\sigma,\lambda)\mapsto \mrm{IC}(\sigma,\lambda)$
defines a map $\Delta_{\nu} \to \mathcal B_V$.  
This map is injective due to the fact that different pairs 
$(\sigma,\lambda)$ determine different perverse sheaves. 
Moreover, the cardinalities
$|\Delta_{\nu}|\overset{(1)}{=}|\mrm{Irr}\;\Lambda_V|\overset{(2)}{=} |\mathcal P_V|$,
where $\mrm{Irr}\;\Lambda_V$ is the set of all irreducible component of the variety
$\Lambda_V$ constructed in \cite{lusztig1}. The equality (1) holds by
\cite[corollary 5.3]{Ringel2} and 
the equality (2) holds by \cite[Theorem 4.16 (b)] {lusztig2}.
By definitions, the two sets are of finite order. 
Therefore we have

\begin{thm}
(\cite[Theorem 6.16 (b)]{lusztig2})
The map $\Delta_{\nu} \to \mathcal B_V$ is bijective.
\end{thm}

\subsection{Proof of Proposition ~\ref{generalcases}}
We preserve the setting of Section ~\ref{general}.
Given any element $(\sigma,\lambda)\in \Delta_{\nu}$,
recall that the variety $X(\sigma, \lambda)$ contains all elements
$x\in E_{V,\Omega}$ such that
\[
(V,x)\simeq \oplus M^{\sigma(M)} \, \oplus\, (R_1 \oplus \cdots \oplus R_q),
\]
where $R_1,\cdots R_q$ are pairwise nonisomorphic homogeneous regular 
simple representations in $\mrm{Ind}(Q)$.

We can write the representation 
$\oplus M^{\sigma(M)} \, \oplus\, (R_1 \oplus \cdots \oplus R_q)$ as
\begin{align*}
&O_1 \oplus \cdots \oplus O_n,
\end{align*}
such that
\begin{align*}
&\mrm{Ext}^1(O_m, O_{m'})=0\;\text{and}\; \mrm{Hom}(O_{m'}, O_m)=0,\;\text{
  if}\; m < m',
\tag{a}
\end{align*}
and $O_m$ has one of the following forms:
\begin{align*}
&O_m=M^{\sigma(M)}\quad \text{ where}\; M \in \text{Ind}(Q) \;
\text{is preprojective or preinjective} ;\tag{1}\\
&O_m=\oplus_{M\in T} \, M^{\sigma(M)},\quad \text{ where }\; T\;
\text{is a tube};\tag{2}\\
&O_m=R_1\oplus \cdots \oplus R_q.\tag{3}
\end{align*}

The  condition (a) can be accomplished by 
putting the preprojective (resp. preinjective) $O_m$'s in case (1)  
in the first (resp. last) part of the 
sequence and putting the $O_m$'s in cases (2) and (3) in the middle part of the sequence,
then adjusting the $O_m$'s in case (1)  such that they satisfy the 
condition (a). (This can be done due to Lemma ~\ref{vanishing}.) 
Note that any  order of the $O_m$'s in case (2) and (3)
already satisfies the condition (a).
 
For each $m$, let $\nu_m=|O_m|$. 
By the definition of $O_m$, this is well-defined. 
Fix a $K$-vector space $V(m)$ such that $|V(m)|=\nu_m$. 
By abuse of notations, denote by $O_m$ the subvariety in 
$E_{V(m),\Omega}$ consisting of all elements $x$ such that
$(V,x)\simeq O_m$.
(Note that $O_m$ in case (3) is nothing but $X(0)$ in
$E_{V(m),\Omega}$ in Section ~\ref{noncyclic}.)

Define the irreducible local system $\mathscr L_m$ on $O_m$ by
\begin{itemize}
 \item $\mathscr L_m=\bar{\mathbb Q}_l$ when $O_m$ is case (1) or (2);
 \item $\mathscr L_m=\mathscr L_{\chi(\lambda)}$ (see Lemma ~\ref{Lemma-2}) when $O_m$ is case (3).
\end{itemize}
Denote by $\mrm{IC}(O_m,\mathscr L_m)$ the intersection complex
on $E_{V(m),\Omega}$ determined by $O_m$ and $\mathscr L_m$.
Then $\mrm{IC}(O_m,\mathscr L_m)$ is in $\mathcal P_{V(m)}$ by Lemma
~\ref{Lemma-1},~\ref{tube} and ~\ref{Lemma-2} for 
$O_m$ in the case (1), (2) and  (3), respectively.
So  the semisimple complex
\[
\mrm{IC}(O_1,\mathscr L_1)\star\cdots\star\mrm{IC}(O_n,\mathscr L_n)
\]
on $E_{V, \Omega}$ is in $\mathcal Q_V$ (see Section ~\ref{inductionfunctor}).

To prove Proposition ~\ref{generalcases}, it suffices to show that 
$\mrm{IC}(\sigma,\lambda)$ is a direct summand of the semisimple complex
$\mrm{IC}(O_1,\mathscr L_1)\star\cdots\star\mrm{IC}(O_n,\mathscr L_n)$
up to shift.

For simplicity, denote by $E_m$ the variety $E_{V(m),\Omega}$ for
$m=1,\cdots, n$. Let $\bar O_m$ the closure of $O_m$ in $E_m$ for
$m=1,\cdots,n$.
Recall from Section ~\ref{inductionfunctor},
we have the following diagram
\[
\begin{CD}
E_1\times \cdots \times E_n @<q_1<< F' @>q_2>> F'' @>q_3>> E_{V,\Omega}.
\end{CD}
\tag{**}
\]
By Lemma ~\ref{n-mult}, we have
\[
\mrm{IC}(O_1,\mathscr L_1)\star\cdots\star\mrm{IC}(O_n,\mathscr L_n)
=(q_3)_! (q_2)_{\flat} (q_1)^*(\mrm{IC}(O_1, \mathscr L_1)\boxtimes \cdots \boxtimes
\mrm{IC}(O_n,\mathscr L_n)).
\tag{4}
\]

Let $\tilde{\mathcal F}''$ be the subvariety of $F''$ consisting of 
all elements $(x, V^{\bullet})$ such that the induced representations
$(V^{m-1}/V^m, x)$ is in $\bar{O}_m$ 
for $m=1,\cdots,n$.

Let $\tilde{\mathcal F}'$ be the subvariety of $F''$ consisting of 
all elements $(x, V^{\bullet})$ such that the induced representations
$(V^{m-1}/V^m,x)$ is in $O_m$ for any $m=1,\cdots,n$.

Denote by $A''$ the subvariety of $F'$ consisting of all triples $(x,V^{\bullet},\mbf g)$ 
in $E'$ such that the induced representations of $x$ are in  the $\bar{O}_m$'s.

Denote by $A'$ the subvariety of $F'$ consisting of all triples $(x, V^{\bullet},\mbf g)$
such that the induced representations of $x$ are in the $O_m$'s.

Consider the following commutative diagram:
\[
\begin{CD}
 O_1\times\cdots\times O_n @<q_1''<<  A'       @>q_2''>> \tilde{\mathcal F}'\\
 @Vj'VV                              @Vi'VV                 @ViVV\\
 \bar{O}_1\times\cdots\times \bar{O}_n  @<q_1'<<  A'' @>q_2'>> \tilde{\mathcal F}''\\
 @VjVV                                           @V\rho 'VV           @V\rho VV\\
E_1\times\cdots\times E_n @<q_1<<   F' @>q_2>>        F'' @>q_3>> E_{V,\Omega} 
\end{CD}
\]
where the bottom row is the diagram (**),
$q_1'$ and $q_1''$ are the restrictions of $q_1$, 
$q_2'$ and $q_2''$ are the restrictions of $q_2$, and 
the vertical maps are natural embeddings.

From the definitions, the squares in the above diagram are Cartesian.
Since $q_3$  and $\pi:=q_3\rho$ are proper, $\rho$ is proper. Hence $\rho'$ is proper.
By the base change theorem for proper morphisms, we have 
\[
q_1^* j_!=\rho_! (q_1')^* \quad \text{and} \quad q_2^* \rho_!=\rho'_!(q_2')^*.
\]
Note that 
\[
\boxtimes_{m=1}^n \mrm{IC}(O_m,\mathscr L_m)=j_{!*} j'_{!*} (\boxtimes
\mathscr L_m)[\dim O_1\times\cdots \times O_n].
\]
Set $d=\dim O_1\times \cdots \times O_n$. Note that $j_{!*}=j_!$.
So
\begin{equation*}
\begin{split}
&(q_3)_! (q_2)_{\flat} (q_1)^* (\boxtimes_{m=1}^n\mrm{IC}(O_m,\mathscr L_m))
=(q_3)_! (q_2)_{\flat} (q_1)^* j_! j_{!*}' (\boxtimes_{m=1}^n
\mathscr L_m)[d]\\
&=(q_3)_! (q_2)_{\flat} \rho_! (q_1')^* (j_{!*}' (\boxtimes_{m=1}^n
\mathscr L_m)[d]
=(q_3)_! \rho_! (q_2')_{\flat} (q_1')^* j_{!*}'(\boxtimes_{m=1}^n
\mathscr L_m)[d]\\
&=\pi_! (q_2')_{\flat} (q_1')^* j_{!*}' (\boxtimes_{m=1}^n \mathscr L_m)[d].
\end{split}
\end{equation*}
Denote by $\mathscr L$ the complex 
$(q_2')_{\flat} (q_1')^* j_{!*}' (\boxtimes_{m=1}^n \mathscr L_m)[d]$.
So (4) becomes
\[
(q_3)_! (q_2)_{\flat} (q_1)^* (\boxtimes_{m=1}^n\mrm{IC}(O_m,\mathscr
L_m))=\pi_!(\mathscr L).
\tag{5}
\]
Note that $j_{!*}(\boxtimes_{m=1}^n\mathscr L_m)[d]$ is a simple
perverse sheaf. Recall that $q_1'$ is smooth with connected fibres, by
\cite[Proposition 4.2.5]{BBD},
\[
(q_1')^* [d_1](j_{!*}(\boxtimes_{m=1}^n(\mathscr L_m))[d] \quad \text{ is
  a simple perverse sheaf on $A''$.}
\]
Since $q_2'$ is a principal $G_{V(1)}\times \cdots \times G_{V(n)}$-bundle, 
$\mathscr L$  is a simple perverse sheaf on $\tilde{\mathcal F}''$ up
to shift.
Note that $O_1\times \cdots \times O_n$ is a smooth variety, by 
\cite[Lemma 4.3.2]{BBD}, 
\[
(j')^* j'_{!*}(\boxtimes_{m=1}^n\mathscr L_m)=\boxtimes_{m=1}^n
\mathscr L_m.
\]
Since the top square in the above diagram are Cartesian, 
\[
i^*\mathscr L=(q_2'')_{\flat} (q_1'')^* (j')^*
j'_{!*}(\boxtimes_{m=1}^n\mathscr L_m)[d]=
(q_2'')_{\flat} (q_1'')^* (\boxtimes_{m=1}^n \mathscr L_m)[d].
\tag{6}
\]
We set
$X=X(\sigma, \lambda)-X(\sigma, \lambda)\cap 
\pi(\tilde{\mathcal F}''-\tilde{\mathcal F}')$ and
$Y=\pi^{-1}(X)$. Then
\[
\text{$X$ is open dense in $X(\sigma,\lambda)$ and the restriction 
$\pi^0:Y\to X$ is an isomorphism.} 
\]
(See [L2, 6.12 (a) (b)] or [Li, 5.6], this is where the condition (a) is used.)
Thus we have the following diagram:
\[
\begin{CD}
O_1\times\cdots\times O_n @<q_1^0<<   F^0  @>q_2^0>>
 Y @>\pi^0>>  X,
\end{CD}
\]
where $F^0=(q_2')^{-1}(Y)$ and $q_1^0$, $q_2^0$, $\pi^0$ are the natural restrictions of 
$q_1'$, $q_2'$ and $\pi$, respectively.
By the definition of $\mathscr L_{\chi(\lambda)}$, we have
\[
(\pi^0 q_2^0)^*(\mathscr L_{\chi(\lambda)}|_{X})
=(q_1^0)^*(\boxtimes_m\mathscr L_m).
\] 
Also by (6), $(q_1^0)^*(\boxtimes_m\mathscr L_m)=(q_2^0)^*\mathscr L|_Y$.
So $(\pi^0)^*(\mathscr L_{\chi(\lambda)}|_X)=\mathscr L|_Y$.
Since $\pi^0:Y\to X$ is an isomorphism, 
\[
(\pi)_!(\mathscr L)[-d]|_X \simeq \mathscr L_{\chi(\lambda)}|_X.\] 
Therefore, 
the intersection complex $\mrm{IC}(X,\mathscr L_{\chi(\lambda)}|_X)$ is a direct summand
of $\pi_!(\mathscr L)[d]$. 
Since $X$ is open dense in the closure of $X(\sigma,\lambda)$,
$\mrm{IC}(X,\mathscr L_{\chi(\lambda)}|_X)=\mrm{IC}(\sigma,\lambda)$.
Proposition ~\ref{generalcases} follows.

\section{Comments}
Note that we deal with  the characterizations of the canonical bases in the symmetric cases.
It may be of interest to characterize the canonical bases in the
nonsymmetric cases. 

From the proof of Proposition ~\ref{generalcases}, the set
\[
\{\mrm{IC}(O_1,\mathscr L_1)\star \cdots \star \mrm{IC}(O_n,\mathscr
L_n)\;|\; (\sigma, \lambda) \in \Delta_{\nu}, \nu\in \mathbb N[I]\}
\]
is a $\mathbb Q(v)$-basis of the algebra $(\mathcal K_Q, \star)$ (see ~\ref{algebra}).
Moreover by looking closer to the shifts, one can show that 
\begin{cor}
The set
\[
\{
\mrm{IC}(O_1,\mathscr L_1)\circ \cdots \circ \mrm{IC}(O_n,\mathscr
L_n)\;|\; (\sigma,\lambda)\in \Delta_{\nu},\nu\in \mathbb N[I]
\}
\]
is an $\mathbb A$-basis of the algebra 
$(\mathcal K,\circ)=\,_{\mathbb A}\mbf U^-$ (see ~\ref{algebra}) and 
stable under bar involution. The
transition matrix between this basis and the canonical basis is upper triangular with
entries in the diagonal equal 1 and entries above the diagonal in
$\mathbb A$. 
\end{cor}

\begin{proof}
For simplicity, we write $\mathbf C_{\sigma,\lambda}$ for the complex
\[
\mrm{IC}(O_1,\mathscr L_1)\circ \cdots \circ \mrm{IC}(O_n,\mathscr L_n).
\]
Since the bar involution commutes with the multplication $\circ$ (see \cite{lusztig3})
and the intesection cohomology complexes are self-dual, 
the complex $\mathbf C_{\sigma, \lambda}$ is stable under the bar involution.
From the proof of Proposition ~\ref{generalcases}, we see that 
\[
\mathbf C_{\sigma,\lambda}=\mrm{IC}(\sigma,\lambda)[d] \oplus P,
\]
for some $d$ and 
$supp(P)\subseteq \overline{X(\sigma,\lambda)}- X(\sigma, \lambda)$.
But $\mathbf C_{\sigma,\lambda}$ is bar invariant, $d$ has to be zero. Corollary follows.
\end{proof}

The relationship between this
basis and the ``canonical basis'' defined in ~\cite{LXZ} deserves
further investigation.

\end{document}